\documentclass[11pt,final]{article} 

\usepackage{amsfonts,amsthm,amsmath, amssymb, mathtools, bm}
\usepackage{graphicx,color}
\usepackage{hyperref}
\usepackage{a4wide}
\usepackage{multirow}
\usepackage{enumitem}
\usepackage{tikz-cd}
\usepackage{subcaption}
\usepackage{mathrsfs}
\usepackage{pgfplots}
\pgfplotsset{compat=1.18}
\usepackage{tikz}

\usepackage{aliascnt}
\usepackage{cleveref}
\newtheorem{proposition}{Proposition}[section]

\newaliascnt{theorem}{proposition}
\newtheorem{theorem}[theorem]{Theorem}
\aliascntresetthe{theorem}

\newaliascnt{lemma}{proposition}
\newtheorem{lemma}[lemma]{Lemma}
\aliascntresetthe{lemma}

\newaliascnt{corollary}{proposition}

\aliascntresetthe{corollary}

\newaliascnt{definition}{proposition}
\newtheorem{definition}[definition]{Definition}
\aliascntresetthe{definition}

\newaliascnt{remark}{proposition}
\newtheorem{remark}[remark]{Remark}
\aliascntresetthe{remark}

\newaliascnt{example}{proposition}
\newtheorem{example}[example]{Example}
\aliascntresetthe{example}

\crefname{definition}{definition}{definitions}
\Crefname{definition}{Definition}{Definitions}

\crefname{theorem}{theorem}{theorems}
\Crefname{theorem}{Theorem}{Theorems}

\crefname{lemma}{lemma}{lemmas}
\Crefname{lemma}{Lemma}{Lemmas}

\crefname{corollary}{corollary}{corollaries}
\Crefname{corollary}{Corollary}{Corollaries}

\crefname{remark}{remark}{remarks}
\Crefname{remark}{Remark}{Remarks}

\crefname{example}{example}{examples}
\Crefname{example}{Example}{Examples}

\crefname{figure}{figure}{figures}
\Crefname{figure}{Figure}{Figures}

\newenvironment{proofof}[1]{\smallskip\noindent{\textbf{Proof~of~#1.}}%
  \hspace{1pt}}{\hspace{-5pt}{\nobreak\quad\nobreak\hfill\nobreak%
    $\square$\vspace{2pt}\par}\smallskip\goodbreak}

\numberwithin{equation}{section}
\numberwithin{figure}{section}

\allowdisplaybreaks[4]

\renewcommand{\phi}{\varphi}
\renewcommand{\theta}{\vartheta}
\renewcommand{\epsilon}{\varepsilon}
\renewcommand{\L}[1]{\mathbf{L^#1}}

\renewcommand{\d}[1]{\mathinner{\mathrm{d}{#1}}}

\renewcommand{\L}[1]{\mathbf{L^#1}}

\newcommand{\Lloc}[1]{\mathbf{L^{#1}_{loc}}}
\newcommand{\Wloc}[2]{\mathbf{W^{#1,#2}_{loc}}}
\newcommand{\C}[1]{\mathbf{C^{#1}}}

\newcommand{\Cc}[1]{\mathbf{C_c^{#1}}}

\newcommand{\W}[2]{\mathbf{W^{#1,#2}}}
\newcommand{\BV}{\mathbf{BV}}
\newcommand{\cBV}[1]{\mathcal{BV}^{#1}}
\newcommand{\BVloc}{\mathbf{BV_{loc}}}
\newcommand{\modulo}[1]{{\left|#1\right|}}
\newcommand{\norma}[1]{{\left\|#1\right\|}}
\newcommand{\reali}{{\mathbb{R}}}
\newcommand{\naturali}{{\mathbb{N}}}
\newcommand{\interi}{{\mathbb{Z}}}
\newcommand{\Lip}{\mathop\mathbf{Lip}}
\newcommand{\tv}{\mathop\mathrm{TV}}

\newcommand{\esssup}{\mathop\mathrm{ess~sup}}

\newcommand{\sgn}{\mathop\mathrm{sgn}}
\newcommand{\spt}{\mathop\mathrm{spt}}
\newcommand{\Id}{\mathop{\mathbf{Id}}}
\newcommand{\pint}[1]{\mathaccent23{#1}}

\newcommand{\caratt}[1]{\chi_{\strut#1}}

\newcommand{\Su}{\mathbb{S}^1}

\DeclareFontFamily{U}{mathx}{\hyphenchar\font45}
\DeclareFontShape{U}{mathx}{m}{n}{
  <-6> mathx5 <6-7> mathx6 <7-8> mathx7
  <8-9> mathx8 <9-10> mathx9
  <10-12> mathx10 <12-> mathx12
}{}
\DeclareSymbolFont{mathx}{U}{mathx}{m}{n}
\DeclareMathAccent{\widecheck}{0}{mathx}{"71}

\begin{document}

\title{Traveling Waves for Mixed Nonlocal-Local \\ Systems of
  Conservation Laws on $\Su$}

\author{Rinaldo M.~Colombo$^1$ \and Mauro Garavello$^2$ \and Claudia
  Nocita$^2$}

\maketitle

\footnotetext[1]{Unità INdAM \& Dipartimento di Ingegneria
  dell'Informazione, Università di Brescia, Italy.\hfill\\
  \texttt{rinaldo.colombo@unibs.it}}

\footnotetext[2]{Dipartimento di Matematica e Applicazioni,
  Università di Milano-Bicocca, Italy.\hfill\\
  \texttt{mauro.garavello@unimib.it} and
  \texttt{c.nocita@campus.unimib.it}}

\begin{abstract}
  \noindent Non trivial traveling waves in genuinely nonlinear
  conservation laws hardly appear and, in some cases, scalar nonlocal
  conservation laws inherit this feature. Here, we show that mixed
  nonlocal-local systems of conservation laws display a far richer
  scenario. A variety of traveling waves may exist, different
  components may travel with different speeds and traveling waves can
  appear as time-asymptotic limits. Technically, these results need a
  significant extension of Dafermos' \emph{``trend to steady state''}
  theory, exploiting properties of Hamiltonian systems solved by
  generalized characteristics.

  \medskip

  \noindent\textit{2020~Mathematics Subject Classification:}
  35L65; 35C07; 35B10; 34C37

  \medskip


\end{abstract}


\section{Introduction}
\label{sec:introduction}

We aim at investigating periodic traveling wave profiles for mixed
systems of nonlocal and local conservation laws in one space dimension
motivated by traffic modeling.  Given $k \in \naturali$, speed
functions $v_{NL}^i, v_L \colon \reali \to \reali$ and convolution
kernels $\eta^i \colon \reali \to \reali^{k+1}$, we consider the
following Cauchy problem:
\begin{equation}
  \label{eq:3}
  \begin{cases}
    \partial_t \rho^i
    +
    \partial_x\left(\rho^i \;
    v_{NL}^i\left(\left(\Sigma\rho + r\right)*\eta^i\right) \right)= 0
     & \text{in} \, \reali_+ \times \Su
    \qquad i \in \{1, \ldots, k\}
    \\
    \partial_t r +
    \partial_x \left(r \;
    v_L\left(\Sigma\rho + r\right)\right) =0
     & \text{in} \, \reali_+ \times \Su
    \\
    \rho(0,x)=\rho_o(x)
     & \text{in} \,  \Su
    \\
    r(0,x)=r_o(x)
     & \text{in} \,  \Su \,.
  \end{cases}
\end{equation}
The space periodic unknown
$(\rho, r): \reali_+ \times \Su \to \reali^{k+1}$ consists of $k$
\emph{nonlocal} components $\rho^i$ and a \emph{local} one $r$.
Nonlocal and local equations in~\eqref{eq:3} are coupled in the speeds
through the total density
\begin{equation}
  \label{eq:60}
  \Sigma\rho + r \coloneqq \sum_{i=1}^k \rho^i + r.
\end{equation}
While the speed $v_L$ is a function of the pointwise value of
$\Sigma \rho +r$, $v_{NL}^i$ depends on its convolution with a kernel
$\eta^i \in \L1(\reali; \reali^{k+1})$ for $i \in \{1, \ldots, k\}$:
\begin{equation}
  \label{eq:98}
  \left(\Sigma \rho + r\right)*\eta^i
  \coloneqq
  \sum_{j=1}^k \rho^j*\eta^i_j + r*\eta^i_{k+1}\,.
\end{equation}
We defer to~\eqref{eq:14} for the precise definition of periodic
convolution.

System~\eqref{eq:3}, introduced in~\cite{opera-seconda} and extended
in~\cite{periodic} to the periodic setting, models traffic composed by
human-driven vehicles interacting with $k$ classes of autonomous
vehicles (AVs): $r$ represents the density of human-driven vehicles,
$\rho^i$ that of the $i$-th class of AVs.  While AVs are aware of
traffic conditions over long distances, human drivers perceive them
only locally.  This difference is captured by the nonlocal-local
coupling in~\eqref{eq:3}.

The phenomenon of \emph{stop and go} waves in traffic flow is widely
studied for its impact on congestion, fuel consumption and driving
safety, see for instance~\cite{Sugiyamal2008, Tadaki2013}. In this
context, AVs have recently been proposed as stabilizers of traffic,
through strategies designed to dissipate such
waves~\cite{STERN2019351, STERN2018205}.  Motivated by these
experiments, typically carried on a ring road, we analytically
investigate traveling wave profiles for~\eqref{eq:3} in a periodic
setting.

Despite their relevance in applications, a limited number of results
regarding traveling waves for nonlocal conservation laws is available
in the literature.  Most concern kernels supported on the negative
semiaxis, so that $v_{NL}$ depends on a weighted average of the
downstream traffic.  The existence of nonconstant, monotone, $\C1$
traveling wave profiles is established in~\cite{MR3960494} for a
scalar nonlocal conservation law on the real line, together with their
property of asymptotic-in-time attractors.  This result is extended
in~\cite{MR4029531, MR4026006} to the case of rough roads and to a
slow erosion model in~\cite{MR3115842}.  For scalar periodic
conservation laws with linear nonlocal speeds and kernels increasing
on the negative semiaxis, \cite{MR4371088} ensures that solutions
$\L2$-converge as $t \to +\infty$ to the mean value of the initial
datum.  An analogous result is proven in~\cite{MR4384589} in the case
of possible nudging terms, dropping the anisotropic assumption on the
kernel.  In applications to traffic management, these results indicate
the possibility of dissipating traffic waves through a suitable choice
of $v_{NL}$ and $\eta$.  Control-oriented results include the exact
boundary controllability of~\cite{MR4298705} and the Lyapunov
stabilization of~\cite{MR4643838}.  See~\cite{MR3013408} for traveling
waves in balance laws with a nonlocal source term. The Camassa-Holm
equation, which can be seen as a further example of nonlocal
conservation law, is treated in~\cite{MR4387318}.

For local systems of conservation laws, classical results --- in
particular decay estimates --- confine the possible appearance of
traveling waves to very particular cases, see for
instance~\cite{MR454375} or the review~\cite{MR580309}.  In
particular, strictly hyperbolic, genuinely nonlinear systems only
admit trivial traveling waves.

We present the first analysis of traveling wave profiles in mixed
nonlocal-local systems, under the mild assumptions
$v_{NL} \in \L\infty(\reali; \reali^k)$,
$v_L \in \Lip(\reali; \reali)$ non increasing and
$\eta \in \L1(\reali; \reali^{k \times (k+1)})$, so that the study
below includes the settings of~\cite{MR4371088, MR4384589, MR3960494}.

In \Cref{prop:traveling} we show that~\eqref{eq:3} admits a large
variety of non trivial (nonconstant) traveling waves evolving with
speed $\lambda \in \reali$, whenever the nonlocal and local speeds
evaluated on the profile equal $\lambda$.  Under the stronger
assumptions of \Cref{thm:WP}, each of these profiles is in fact the
unique solution to~\eqref{eq:3}.  Also, these profiles are not
asymptotically stable. Thus, nonlocality makes the traveling waves
scenario dramatically richer than that related to local systems.

The mere assumption that $v_L$ be weakly monotone does not prevent
unexpected behaviors of traveling waves, so that the components of
$\rho$ and $r$ may evolve with different speeds, see \Cref{ex:tw}.
The stronger strict monotonicity assumption of $v_L$ rules out this
occurrence, see \Cref{lem:Mauro} and \Cref{thm:waves}.

Assuming that $r_o$ is almost everywhere non zero and that
$v_L^{-1}(\lambda)$ has empty interior, the presence of a solution
whose local component $r$ is a traveling wave with speed $\lambda$
forces $\Sigma \rho$ to be a traveling wave with the same speed
$\lambda$, see \Cref{thm:waves}.

On the other hand, the converse question, i.e., whether $r$ inherits
the traveling wave structure of $\Sigma \rho$, admits a weaker,
asymptotic-in-time answer. Assuming that the restriction of $v_L$ to
some interval $[0,R]$ is of class $\C2$, strictly decreasing and
concave, \Cref{thm:Claudia} shows that whenever $\rho$ is a traveling
wave with speed $\lambda$, $r$ converges almost uniformly as
$t \to + \infty$ to a traveling wave $r_\infty$ which we describe in
detail. Indeed, $r_\infty$ travels with the same speed $\lambda$ and
displays finitely many jump discontinuities connecting two continuous
profiles explicitly characterized and depending only on
$\Sigma \rho_o$, $v_L$ and $\lambda$.

Inspired by~\cite{Trend-to-steady}, we consider the Hamiltonian system
solved by the generalized characteristics of the local equation
in~\eqref{eq:3} and use heteroclinic and homoclinic curves to localize
the long-time behavior of $r$. However, the structure of the present
Hamiltonian system is significantly more intricate than that
considered in~\cite{Trend-to-steady}. Indeed, we here face a phase
portrait that may feature any finite number of saddles and centers,
while~\cite{Trend-to-steady} deals with exactly one center and two
saddles. Due to this rich structure, a rather careful analysis is
necessary to show that, as $t$ grows, the shocks in $r (t)$ either
annihilate or converge to the jumps in $r_\infty (t)$.

Remarkably, the core of \Cref{thm:Claudia} is actually a result on the
asymptotic behavior of solutions to a scalar (local) conservation law
with $x$-dependent flux, obtained extending the techniques
in~\cite{Trend-to-steady} as detailed above. Indeed, the role of the
mixed nonlocal-local structure of~\eqref{eq:3} simply enters the
general dependence of the scalar flux on the space variable.

\medskip

The next Section presents the analytical results. All proofs are
collected in Section~\ref{sec:proofs}, while general technical results
are deferred to the Appendix.

\section{Main Results}\label{sec:main-result}
We represent elements of $\Su$ with real numbers in
$\mathopen[0,1\mathclose[$. Given a map $f \colon \Su \to \reali$ we
denote by $\widehat f \colon \reali \to \reali$ its $1$-periodic
extension to $\reali$. Conversely, given a $1$-periodic function
$g \colon \reali \to \reali$ by $\widecheck{g} \colon \Su \to \reali$
we denote its restriction to $\Su$. We defer the rigorous definitions
of function spaces, norms, and of the convolution in the periodic
setting to Section~\ref{sec:proofs}.

The space $\cBV{}(\Su; \reali^m)$, for $m \in \naturali$, denotes the
space of functions with bounded total variation on $\Su$ and
$\cBV{1}(\Su; \reali^m)$ is the set of functions in
$\W11(\Su; \reali^m)$ whose spatial derivative has bounded total
variation, see~\eqref{eq:20}. We understand that for any function in
$\cBV{} (\Su; \reali^m)$ we use its right continuous representative.

\medskip

We recall the definition of solutions to~\eqref{eq:3}
presented in~\cite[Definition~2.1]{periodic}.

\begin{definition}
  \label{def:sol}
  Let $I = \reali_+$ or $I=[0,T]$ for a positive $T$. The couple
  \begin{displaymath}
    (\rho,r) \in
    \C0\left(I; \W11 (\Su; \reali^k) \times \L1 (\Su; \reali)\right)
  \end{displaymath}
  is a \emph{solution} to~\eqref{eq:3} on $I$ if, setting for
  $(t,x,u) \in I \times \reali \times \reali$ and for
  $i \in \{1, \ldots, k\}$
  \begin{eqnarray*}
    v^i (t,x)
    & \coloneqq
    & v_{NL}^i
    \left(\left(\Sigma\widehat{\rho (t)} * \eta^i\right) (x)
    +
    \left(\widehat{r(t)} *\eta^{i}_{k+1}\right) (x)\right)
    \\
    w (t,x,u)
    & \coloneqq
    & v_L\left(\Sigma\widehat{\rho (t)} \, (x) + u \right),
  \end{eqnarray*}
  then for all $\phi \in \Cc1(\pint{I} \times \reali; \reali)$,
  $\psi \in \Cc1(\pint{I} \times \reali; \reali_+)$,
  $\kappa \in \reali$ and $i \in \{1, \ldots, k\}$
  \begin{eqnarray*}
    \int_I \int_{\reali}
    \left(\widehat\rho^i(t,x) \,\partial_t \phi(t,x) + \widehat\rho^i(t,x) \,
    v^i(t,x) \, \partial_x \phi(t,x) \right) \d{x} \d{t}
    & =
    & 0 \,;
    \\
    \int_I \int_{\reali} \modulo{\widehat{r}(t,x) - \kappa} \, \partial_t \psi(t,x)
    \d{x}\d{t}
    \\
    +
    \int_I \int_{\reali} \sgn(\widehat{r}(t,x) - \kappa) \,
    \left[
      \widehat{r}(t,x) \, w\left(t,x,\widehat{r} (t,x)\right)
      - \kappa \, w(t,x,\kappa)
      \right] \partial_x \psi(t,x) \d{x}\d{t}
    \\
    -
    \int_I \int_{\reali} \sgn(\widehat{r}(t,x) - \kappa) \,
    \kappa \, \partial_x w(t,x,\kappa) \, \psi(t,x) \d{x}\d{t}
    & \geq
    & 0 \,;
    \\
    \rho (0)
    & =
    & \rho_o  \,;
    \\
    r (0)
    & =
    & r_o  \,.
  \end{eqnarray*}
\end{definition}

\noindent The minimal hypotheses
$v_{NL} \in \L\infty(\reali; \reali^k)$,
$v_L \in \Lip(\reali; \reali)$ and
$\eta \in \L1(\reali; \reali^{k \times (k+1)})$ ensure that all the
integrals above are finite.

The well-posedness of~\eqref{eq:3} has been established
in~\cite[Theorem~2.2]{periodic} under suitable regularity assumptions on the velocity functions $v_{NL}$ and $v_L$ and
on the convolution kernels $\eta^i$.
For completeness, we report here the precise statement.

\begin{theorem}[Well-Posedness
  of~\eqref{eq:3},{\cite[Theorem~2.2]{periodic}}]
  \label{thm:WP}
  Assume:
  \begin{enumerate}[label=\bf(vNL)]
  \item\label{item:17}
    $v_{NL} \in (\C3\cap \W3\infty)(\reali;\reali^k)$.
  \end{enumerate}
  \begin{enumerate}[label=\bf(vL)]
  \item\label{item:21} $v_L \in \C2(\reali; [0, V_L])$ and there exist
    $V_L, R_L > 0$ such that $v_L(r) = V_L$ for $r \leq 0$ and
    $v_L(r) = 0$ for $r \geq R_L$.
  \end{enumerate}
  \begin{enumerate}[label=$\bm{(\eta)}$]
  \item\label{item:hyp-eta}
    $\eta \in \C3(\reali;\reali^{k\times (k+1)})$ with
    $\mathop{\rm diam} (\spt \eta) \leq H$ for a suitable
    $H \in \naturali \setminus \{0\}$.
  \end{enumerate}
  \noindent Then, the problem~\eqref{eq:3} generates a unique map
  \begin{equation}
    \label{eq:79}
    \mathcal{S} \colon \reali_+ \times
    \cBV{1}(\Su; \reali_+^k)
    \times
    \cBV{}(\Su; [0,R_L])
    \to
    \cBV{1}(\Su; \reali_+^k)
    \times
    \cBV{}(\Su; [0,R_L])
  \end{equation}
  such that
  \begin{enumerate}[label=\textup{\textbf{($\mathcal{S}$\arabic*)}},
    ref=\textup{\textbf{($\mathbf{\mathcal{S}}$\arabic*)}}]
  \item \label{item:29} $\mathcal{S}$ is a global semigroup, i.e.,
    $\mathcal{S}_0 = \Id$ and for all $t,t_* \in \reali_+$,
    $\mathcal{S}_{t} \circ \mathcal{S}_{t_*} = \mathcal{S}_{t+t_*}$.

  \item \label{item:40} The semigroup $\mathcal{S}$ is locally
    Lipschitz continuous with respect to the
    $\W11(\Su; \reali^k) \times \L1(\Su; \reali)$ norm,
    see~\eqref{eq:20}--\eqref{eq:32bis}.

  \item \label{item:23} For any
    $(\rho_o,r_o) \in \cBV{1}(\Su; \reali_+^k) \times \cBV{}(\Su;
    [0,R_L])$ and any $T > 0$, the orbit
    $t \mapsto \mathcal{S}_t (\rho_o, r_o)$ solves~\eqref{eq:3} on
    $[0, T]$ in the sense of \Cref{def:sol}.

  \item \label{item:39} $\mathcal{S}$ is uniquely characterized
    by~\ref{item:29}, \ref{item:40} and by the property that for all
    $(\rho_o, r_o) \in \C2 (\Su; \reali_+^k) \times \cBV{} (\Su;
    [0,R_L])$, the orbit $t \mapsto \mathcal{S}_t (\rho_o,r_o)$ is the
    unique solution to~\eqref{eq:3} in the sense of \Cref{def:sol}.

  \end{enumerate}
\end{theorem}

\noindent We concentrate on traveling wave dynamics
for~\eqref{eq:3}. For $\lambda \in \reali$, by a
\emph{$\lambda$-traveling wave on $\Su$} we mean a function
$f \colon \reali_+ \times \Su \to \reali^m$ ($m \ge 1$) such that, for
a.e. $(t,x) \in \reali_+ \times \reali$,
$\widehat{f}(t,x) = F(x-\lambda \, t) $ for a suitable $1$-periodic
$F \colon \reali \to \reali^m$.  A traveling wave $f$ is \emph{non
  trivial} whenever $F$ is not constant. Note that if
$F \in \Lloc1(\reali; \reali^m)$, then such traveling wave solves, in
distributional sense, the transport equation
\begin{displaymath}
  \partial_t \widehat{f} + \lambda \, \partial_x \widehat{f} = 0
  \quad
  \mbox{in}
  \quad
  \reali_+ \times \reali \, .
\end{displaymath}

In each statement below, we precisely state the assumptions on the
speeds and the kernel functions, typically weaker than those required
in Theorem~\ref{thm:WP}. As a first step, we show that~\eqref{eq:3}
admits a variety of non trivial (not constant) traveling waves.

\begin{proposition}
  \label{prop:traveling}
  Fix $k \!\in \!\naturali \setminus\{0\}$,
  $v_{NL} \in \L\infty(\reali; \reali^k)$,
  $v_L \in \Lip(\reali; \reali)$ and
  $\eta \!\in \! \L1 (\reali; \reali^{k\times (k+1)})$.  
  Fix the initial datum
  $(\rho_o, r_o) \in \W11 (\Su; \reali^k) \times \L1 (\Su; \reali)$
  and assume there exists $\lambda \in \reali$ such that, with the
  notation~\eqref{eq:98},
  \begin{equation}
    \label{eq:1}
    \begin{array}{@{}l}
      \mbox{for all } i \in \{1, \ldots, k\}
      \\
      \mbox{for a.e. } x \in \Su
    \end{array}
    \quad
    \quad
    v_{NL}^i
    \left( \left( \left(\Sigma \rho_o + r_o\right)*\eta^i \right) (x) \right)
    =
    v_L \left( \Sigma\rho_o(x) + r_o (x) \right)
    =
    \lambda\,.
  \end{equation}
  Then, the traveling wave
  $(\rho,r) \colon \reali_+ \times \Su \to \reali^{k+1}$, defined by
  \begin{equation}
    \label{eq:2}
    \left\{
      \begin{array}{rcl}
        \widehat{\rho} (t,x)
        & \coloneqq
        & \widehat{\rho}_o (x-\lambda \, t)
        \\
        \widehat r (t,x)
        & \coloneqq
        & \widehat r_o (x-\lambda \, t)
      \end{array}
    \right.
  \end{equation}
  solves~\eqref{eq:3} in the sense of Definition~\ref{def:sol}.
\end{proposition}

\begin{remark}
  \rm Proposition~\ref{prop:traveling} also provides a wide set of
  \emph{non trivial} traveling wave solutions even in the more
  restrictive case where well-posedness holds by
  Theorem~\ref{thm:WP}. Indeed, under assumptions~\ref{item:17},
  \ref{item:21} and~\ref{item:hyp-eta}, assume that
  $\eta^i_j = \eta_*$ for a fixed $\eta_*$ and
  $i \in \{1, \ldots, k\}$, $j \in \{1, \ldots, k+1\}$. Choose any
  $(\rho_o, r_o) \in \C2(\Su; \reali^k_+) \times \cBV{}(\Su; [0,R_L])$
  such that $\Sigma\rho_o (x) + r_o (x) = \bar s$ for a suitable
  $\bar s \in \reali$. Require that~\eqref{eq:1}, which in the
  present case amounts to
  \begin{displaymath}
    v_{NL}\left( \bar s \int_{\reali} \eta_* (x)\d{x}  \right)
    =
    v_L (\bar s)
    =
    \lambda\,,
  \end{displaymath}
  holds. Then, the unique solution to~\eqref{eq:3} is~\eqref{eq:2},
  which is a non trivial $\lambda$-traveling wave whenever $r_o$ or
  $\rho_o^i$ is not constant.
\end{remark}

\begin{remark}
  \label{rem:as}
  \rm %
  In general, the traveling wave solutions exhibited
  by~\Cref{prop:traveling} are not asymptotically stable. Indeed, fix
  an initial datum $(\rho_o, r_o)$ satisfying~\eqref{eq:1} and
  yielding the traveling wave solution~\eqref{eq:2}. Consider
  arbitrary (non trivial) perturbations
  $(\delta \rho_o, \delta r_o) \in \W11 (\Su; \reali^k) \times \L1
  (\Su; \reali)$ such that $\Sigma \, \delta \rho_o + \delta r_o= 0$
  and
  $\Sigma \, \modulo{\delta \rho_o} + \modulo{\delta r_o} = \epsilon >
  0$.  Then, the traveling wave $(\rho_*, r_*)$ where
  \begin{displaymath}
    \widehat\rho_* (t,x)
    \coloneqq
    \widehat{\rho}_o (x-\lambda\, t)  + \widehat{\delta\rho}_o (x-\lambda\, t)
    \qquad \mbox{ and } \qquad
    \widehat{r}_* (t,x)
    \coloneqq
    \widehat r_o (x-\lambda\, t)  +
    \widehat{\delta r}_o (x-\lambda\, t)
  \end{displaymath}
  solves~(\ref{eq:3}) with the perturbed initial datum
  $(\rho_o + \delta\rho_o, r_o +\delta r_o)$.  $(\rho_*, r_*)$ remains
  distant from the traveling wave solution~\eqref{eq:2} uniformly in
  time, since
  \begin{displaymath}
    \norma{(\rho-\rho_*,r-r_*) (t)}_{\W11(\Su; \reali^k) \times
      \L1(\Su; \reali)} \geq \epsilon
    \qquad \mbox{for all } t \in \reali_+.
  \end{displaymath}
\end{remark}

\smallskip

We begin our study of the analytical properties of traveling wave
solutions to~\eqref{eq:3}. If $\rho$ and $r$ are traveling wave and
$(\rho,r)$ solves~\eqref{eq:3}, $\rho$ and $r$ necessarily evolve with
the same speed.

\begin{proposition}
  \label{lem:Mauro}
  Fix $k \in \naturali \setminus \{0\}$, $R > 0$,
  $v_{NL} \in \L\infty (\reali; \reali^k)$,
  $v_L \in \C2 (\reali; \reali)$ and
  $\eta \in \L1 (\reali; \reali^{k\times (k+1)})$. Assume
  \begin{equation}
    \label{eq:26}
    v_L' (r) < 0 \mbox{ for all } r \in \mathopen]0, R\mathclose[\,.
  \end{equation}
  Let
  $(\rho_o, r_o ) \in \W11 (\Su; \reali_+^k) \times \C0 (\Su;
    \reali_+)$ satisfy
  \begin{equation}
    \label{eq:21}
    \Sigma\rho_o (x) + r_o (x) \in [0, R]
    \quad \mbox{ for all } x \in \Su
    \quad \mbox{ and } \quad
    v''_L \left(\int_{\Su} \left(\Sigma\rho_o (x) + r_o (x)\right)\d{x}\right)
    \leq 0
  \end{equation}
  and assume that~\eqref{eq:3} admits, in the sense of
  Definition~\ref{def:sol}, the solution $(\rho,r)$ where
  \begin{equation}
    \label{eq:22}
    \left\{
    \begin{array}{rcl}
      \widehat{\rho} (t,x)
       & \coloneqq
       & \widehat{\rho}_o (x-\gamma\, t)
      \\
      \widehat r (t,x)
       & \coloneqq
       & \widehat{r}_o (x-c\,t)
    \end{array}
    \right.
    \quad \mbox{ for suitable }\quad
    c,\gamma \in \reali \,.
  \end{equation}
  If $r_o$ is not constant, then $c = \gamma$.
\end{proposition}

Condition~\eqref{eq:26} is sharp: if $v_L$ is non increasing, then the
solution to~\eqref{eq:3} may well consist in a $\lambda$-traveling
wave $\rho$ and a $c$-traveling wave $r$ with $\lambda \neq c$, as
shown in the following example.
\begin{example}
  \label{ex:tw}
  \rm Assume $k=1$, $v_{NL} \in \C0(\reali; \reali)$,
  $v_L \in \Lip (\reali; \reali)$ and $\eta^1_j = \caratt{[0,1]}$ for
  $j \in \{1, 2\}$. Assume that there exist positive $c$, $a$
  and $\Delta$ such that
  \begin{equation}
    \label{eq:12}
    \mbox{for all } r \in [a-\Delta, a+\Delta]
    \qquad
    v_L (r) = c\,.
  \end{equation}
  Then we can choose an initial datum
  $(\rho_o, r_o) \in \W11 (\Su; \reali^k) \times \cBV{} (\Su; \reali)$
  such that~\eqref{eq:3} admits the traveling wave solution $(\rho,r)$
  where, setting $\lambda =v_{NL} (a)$,
  \begin{equation}
    \label{eq:13}
    \left\{
    \begin{array}{rcl}
      \widehat\rho (t,x)
       & =
       & \widehat\rho_o (x - \lambda\,t)
      \\
      \widehat r (t,x)
       & =
       & \widehat r_o (x-c\, t) \,.
    \end{array}
    \right.
  \end{equation}
  In particular, $c$ and $\lambda$ need not be equal and
  $\rho_o + r_o$ is not constant.

  Indeed, choose
  \begin{displaymath}
    \left\{
      \begin{array}{rcl}
        \rho_o (x)
        & =
        & \frac{1}{2} \left( a + \Delta \sin (2\pi\,x) \right)
        \\[4pt]
        r_o (x)
        & =
        & \frac{1}{2} \left( a + \Delta \sin (2\pi\,x) \right)
      \end{array}
    \right.
  \end{displaymath}
  Define $(\rho,r)$ by means of~\eqref{eq:13}. Simple computations
  yield, for all $x \in \reali$,
  \begin{eqnarray*}
    v_{NL}
    \left(
    \left(
    \left( \Sigma\widehat\rho (t) + \widehat r (t) \right) * \eta
    \right)
    (x)
    \right)
    & =
    & v_{NL} (a)
    \\
    v_L\left( \Sigma\widehat\rho (t,x) + \widehat r (t,x) \right)
    & =
    & c
  \end{eqnarray*}
  hence $(\rho,r)$ is a weak solution to~\eqref{eq:3}. Moreover, since
  $\rho$ and $r$ are smooth, also the {Kru\v zkov} entropy inequality
  in Definition~\ref{def:sol} holds. Clearly, $\rho$ is a
  $\lambda$-traveling wave, while $r$ is a $c$-traveling wave.
\end{example}

We now extend the previous proposition to less regular $v_L$ and
$r_o$. The cost is the non vanishing condition~\eqref{eq:10}, which is
however necessary, see~\cref{rem:ro0}.

\begin{theorem}
  \label{thm:waves}
  Fix $k \in \naturali \setminus \{0\}$,
  $v_{NL} \in \L\infty (\reali; \reali^k)$,
  $v_L \in \Lip(\reali; \reali)$ and
  $\eta \in \L1 (\reali; \reali^{k\times (k+1)})$.  Let
  $(\rho_o, r_o) \in \W11 (\Su; \reali^k) \times \cBV{} (\Su; \reali)$
  with
  \begin{equation}
    \label{eq:10}
    r_o (x) \neq 0 \quad \mbox{for a.e. }x \in \Su
  \end{equation}
  yield the map $(\rho,r)$ as solution to~\eqref{eq:3} in the sense of
  \Cref{def:sol}. Assume $r$ is a $\lambda$-traveling wave, for a
  $\lambda \in \reali$, such that
  \begin{gather}
    \label{eq:32}
    v_L^{-1} (\lambda) \mbox{ has empty interior, and}
    \\
    \label{eq:9}
    \lambda \in v_L \left((\Sigma\rho_o+r_o) (\Su)\right) \,.
  \end{gather}
  Then, $\Sigma \rho$ is a $\lambda$-traveling wave and
  $v_L\left(\Sigma \rho(t,x) + r(t,x)\right)= \lambda $ for all
  $(t,x) \in \reali_+\times \Su$.
\end{theorem}

Thus, if the local component $r$ of the solution to~\eqref{eq:3} is
a $\lambda$-traveling wave in $\cBV{}$, then also the sum of the nonlocal
components $\Sigma \rho$ is a $\lambda$-traveling wave.

\begin{remark}
  \label{rm:waves}
  \rm Note that, in the assumptions of~\Cref{thm:waves}, if
  $v_L^{-1}(\lambda)$ is a singleton, then
  \begin{equation}
    \label{eq:65}
    \Sigma \rho(t,x) + r(t,x)= v_L^{-1}( \lambda )
    \quad \mbox{for all } \, (t,x) \in \reali_+\times \Su \, .
  \end{equation}
\end{remark}

Our main result now shows that as soon as the nonlocal component
$\rho$ of the solution to~\eqref{eq:3} is a traveling wave, the local
component $r$ of the solution converges as $t \to +\infty$ to a
traveling. the asymptotic profile of $r$ is fully characterized.

\begin{theorem}
  \label{thm:Claudia}
  Fix $k \in \naturali \setminus \{0\}$,
  $v_{NL} \in \C0 (\reali; \reali^k)$, $v_L \in \Lip(\reali; \reali)$
  and $\eta \in \L1 (\reali; \reali^{k\times (k+1)})$.  Assume there
  exists $R>0$ such that
  \begin{equation}
    \label{eq:26bis}
    {v_L \in \C2([0, R]; \reali)},
    \quad
    v_L' (r) < 0,
    \quad
    v_L'' (r) \leq 0
    \quad
    \mbox{ for all } r \in\mathopen]0, R \mathclose] \,.
  \end{equation}
  Fix
  $(\rho_o, r_o) \in \C2(\Su; \reali^k_+) \times \cBV{}(\Su; \reali)$
  such that
  \begin{gather}
    \label{eq:67}
    \max_{\Su} \Sigma \rho_o \leq R \,;
    \\
    \label{eq:5}
    \left\{ x \in \Su \colon \Sigma\rho_o' (x) = 0 \right\}
    \mbox{ is a finite union of
      closed intervals.}
  \end{gather}
  Assume that~\eqref{eq:3} admits a global solution $(\rho,r)$ in the
  sense of Definition~\ref{def:sol} and that
  \begin{equation}
    \label{eq:68}
    \rho \mbox{ is a } v_L (R)\mbox{-traveling wave.}
  \end{equation}
  Then, there exist
  $\underline{q}_R, \overline{q}_R \in \C0(\Su; \reali)$, uniquely
  characterized below, such that if
  \begin{equation}
    \label{eq:138}
    \underline{q}_R \leq r_o \leq \overline{q}_R\,,
  \end{equation}
  then, $r$ converges asymptotically-in-time almost uniformly in space
  to a $v_L(R)$-traveling wave
  $r_\infty(t,x) \coloneqq q_\infty\left(x-v_L(R)\,t\right)$ where
  \begin{equation}
    \label{eq:33}
    q_\infty(\xi) \coloneqq
    \begin{cases}
      \underline{q}_R(\xi)
      & \xi \in A,
      \\
      \overline{q}_R(\xi)
      & \xi \in \Su \setminus A
    \end{cases}
  \end{equation}
  for suitable $A \subseteq \Su$, defined below. Moreover, $q_\infty$
  is piecewise continuous with finitely many discontinuities, which
  are all positive jumps. Furthermore, $(\rho, r_\infty)$
  solves~\eqref{eq:3} with initial datum $(\rho_o, q_\infty)$ in the
  sense of \Cref{def:sol}.

  The functions $\underline{q}_R$, $\overline{q}_R$ and $q_\infty$ are
  characterized by:
  \begin{enumerate}[label=\bf(\roman*)]

  \item\label{item:43} The functions $\underline{q}_R, \overline{q}_R$
    are uniquely defined by
    \begin{displaymath}
      \mbox{For all } x \in \Su \qquad
      \begin{array}{c}
        0
        \leq \underline{q}_R (x)
        \leq \overline{q}_R (x)
        \leq R - \Sigma \rho_o (x)
        \\
        \underline{q}_R(x) \left(
        v_L\left(
        \underline{q}_R(x) + \Sigma \rho_o(x)
        \right)
        - v_L (R)
        \right)
        = f_R
        \\
        \overline{q}_R(x) \left(
        v_L\left(
        \overline{q}_R(x) + \Sigma \rho_o(x)
        \vphantom{\underline{q}_R}
        \right)
        - v_L (R)
        \right)
        = f_R
      \end{array}
    \end{displaymath}
    where the constant $f_R \in \reali$ is given by
    $ f_R \coloneqq q_R \left( v_L\left(q_R + \max\Sigma\rho_o\right)
      - v_L (R)\right) $ and $q_R \in \reali$ is the unique solution
    to
    $v_L\left(q_R + \max\Sigma\rho_o\right) - v_L (R) + q_R \,
    v_L'\left(q_R {+} \max\Sigma\rho_o\right) =0$.

  \item\label{item:44} To define $q_\infty$, up to a rotation of the
    coordinate on $\Su$, introduce
    $0 = a_0 < \ldots < a_{n} \leq 1$ and
    $0 \leq b_0 < \ldots < b_{n}=1$ with $a_i \leq b_i < a_{i+1}$ so
    that we write~\eqref{eq:5} as
    \begin{equation}
      \label{eq:36}
      \mathop{\rm argmax}_{[0,1]} \Sigma \widehat{\rho}_o
      =
      \bigcup_{i=0}^{n} [a_i, b_i]
    \end{equation}
    and define
    $\left\{ \bar \xi_0, \ldots, \bar \xi_{n-1} \right\} \subset
    [0,1]$ by
    \begin{displaymath}
      \bar \xi_i \in [b_i, a_{i+1} ]
      \quad \mbox{ such that }\quad
      \int_{b_i}^{\bar \xi_i} \underline{q}_R(\xi) \d\xi
      +
      \int_{\bar \xi_i}^{a_{i+1}} \overline{q}_R(\xi) \d\xi
      =
      \int_{b_i}^{a_{i+1}} r_o(\xi) \d\xi.
    \end{displaymath}
    The explicit expression of $q_\infty \colon \Su \to \reali$ is
    \begin{equation}
      \label{eq:33bis}
      q_\infty(\xi) \coloneqq
      \begin{cases}
        \underline{q}_R(\xi)
        & \xi \in \bigcup_{i=0}^{n}[a_i, b_i \mathopen[,
        \\
        \underline{q}_R(\xi)
        & \xi \in \bigcup_{i=0}^{n-1}[b_i, \bar \xi_i \mathopen[,
        \\
        \overline{q}_R(\xi)
        & \xi \in \bigcup_{i=0}^{n-1}[\bar \xi_i, a_{i+1} \mathopen[
      \end{cases}
    \end{equation}
    so that
    $A = \left(\bigcup_{i=0}^{n}[a_i, b_i \mathopen[\right)\cup
    \left(\bigcup_{i=0}^{n-1}[b_i, \bar \xi_i \mathopen[\right)$.
  \end{enumerate}
\end{theorem}

\noindent By \emph{convergence asymptotically-in-time almost uniformly
  in space} we mean that for all sequences $t_\nu$ with
$t_\nu \underset{\nu \to + \infty}{\longrightarrow} +\infty$, for all
$\epsilon>0$ sufficiently small
\begin{displaymath}
  \lim_{\nu \to +\infty} \quad
  \sup_{\{x \in \Su \colon \min_i\modulo{x-v_L(R)t_\nu-\bar \xi_i} > \epsilon\}}
  \modulo{r (t_\nu,x) - q_\infty (x-v_L (R)\, t_\nu)} = 0 \,.
\end{displaymath}

\begin{figure}[!ht]
  \centering
  \includegraphics[width=0.333\textwidth,keepaspectratio,trim=25 25 25
  25]{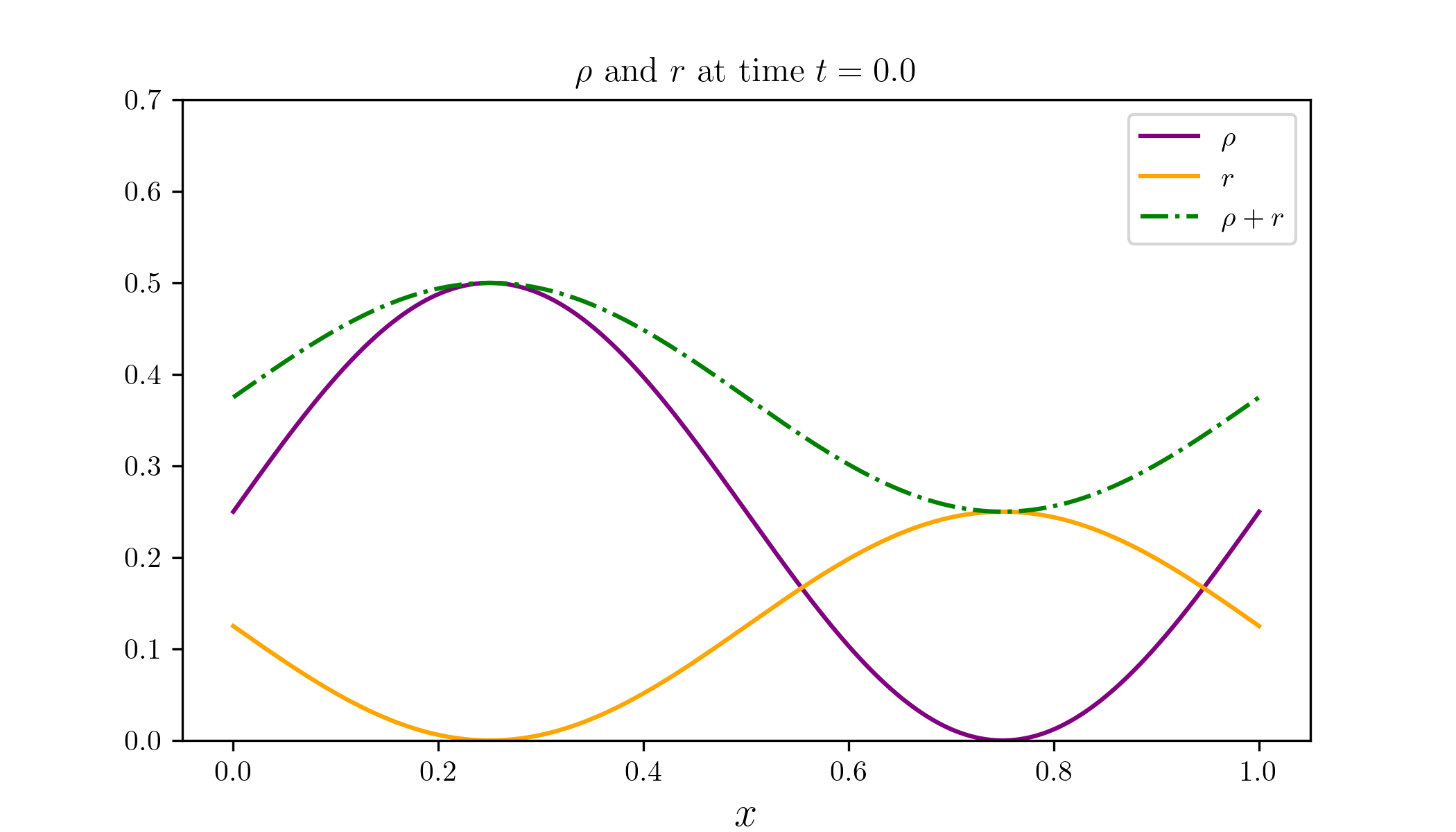}%
  \includegraphics[width=0.333\textwidth,keepaspectratio,trim=25 25 25
  25]{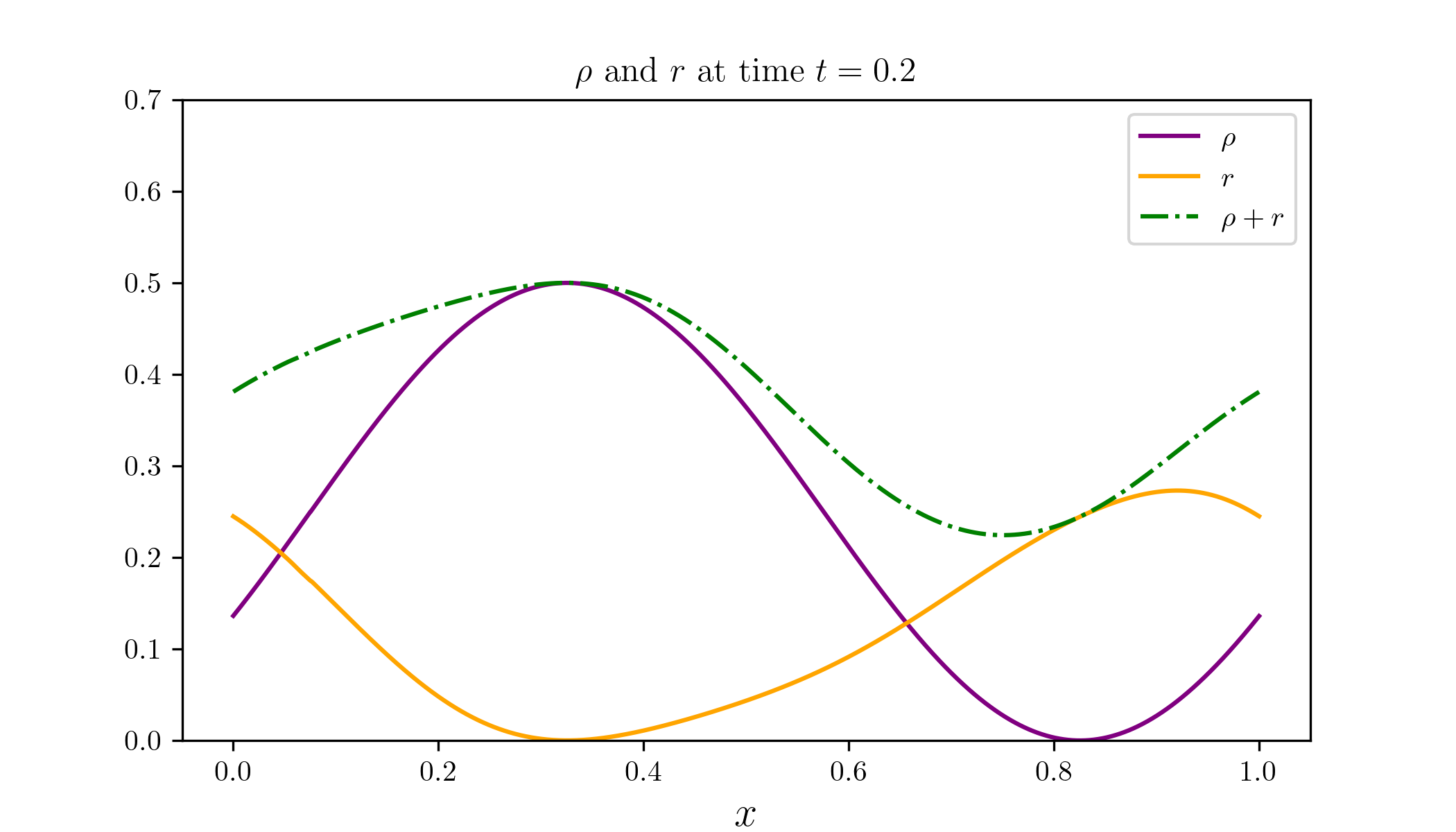}%
  \includegraphics[width=0.333\textwidth,keepaspectratio,trim=25 25 25
  25]{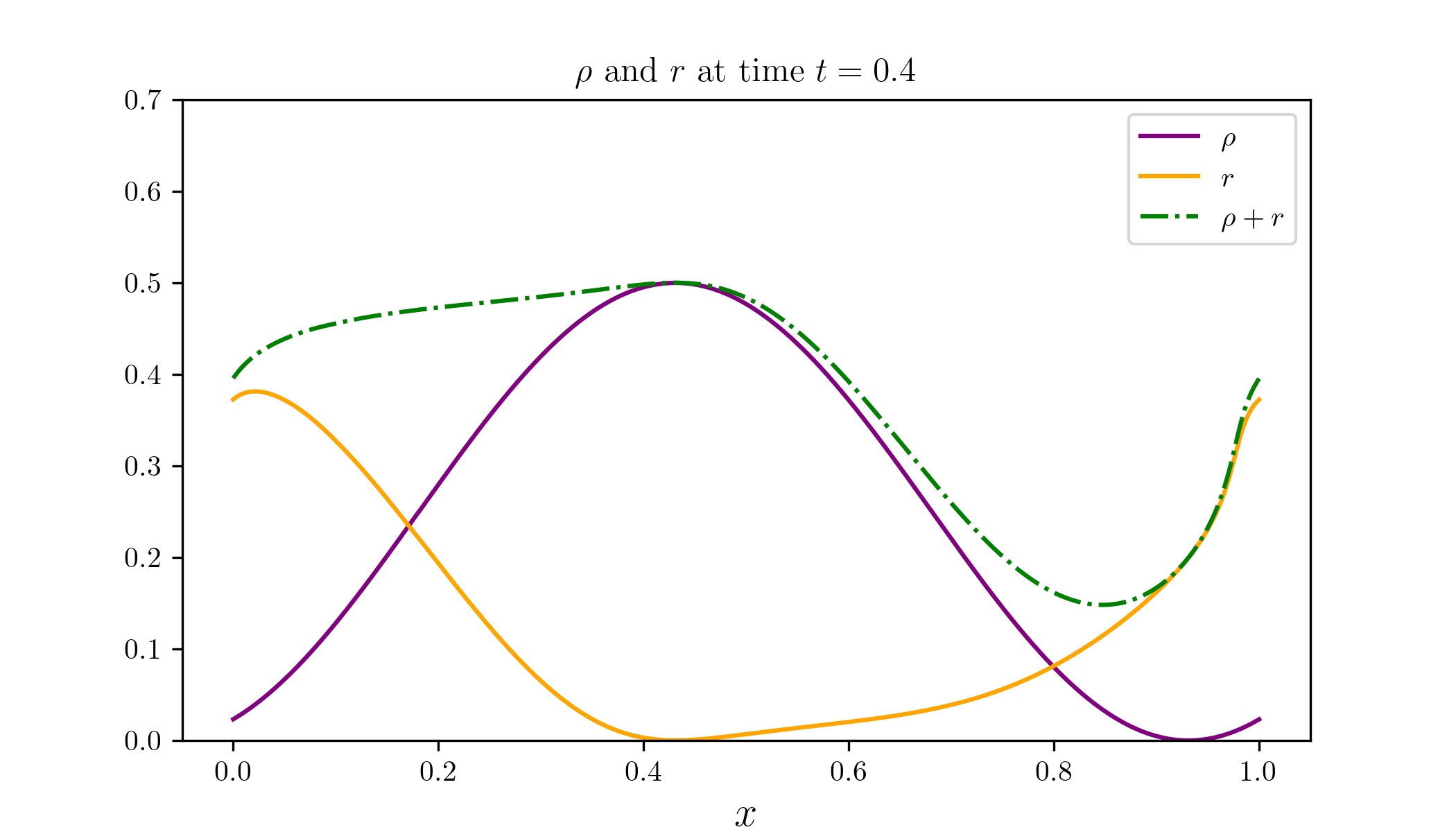}%
  \\[18pt]
  \centering
  \includegraphics[width=0.333\textwidth,keepaspectratio,trim=25 25 25
  25]{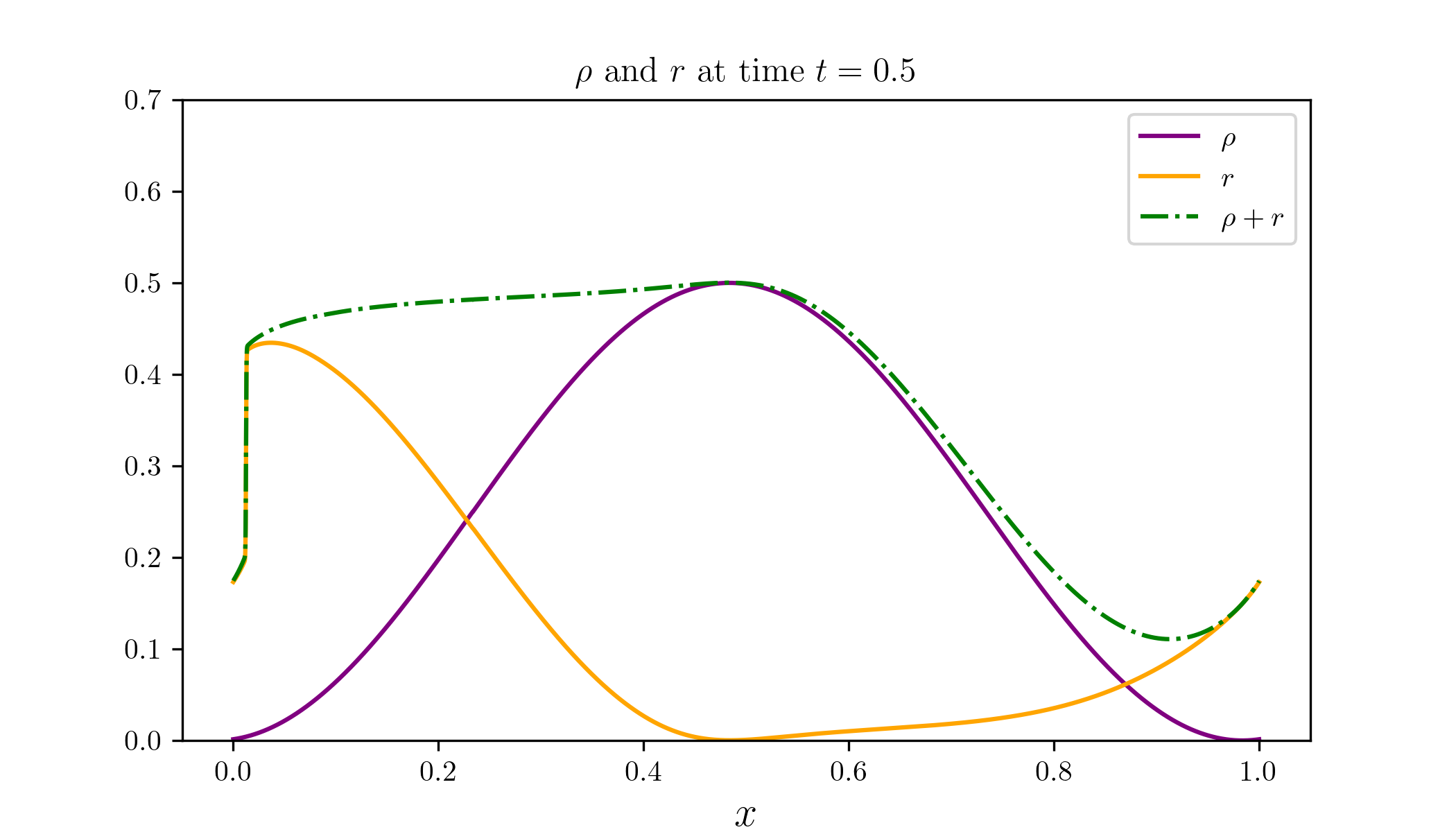}%
  \includegraphics[width=0.333\textwidth,keepaspectratio,trim=25 25 25
  25]{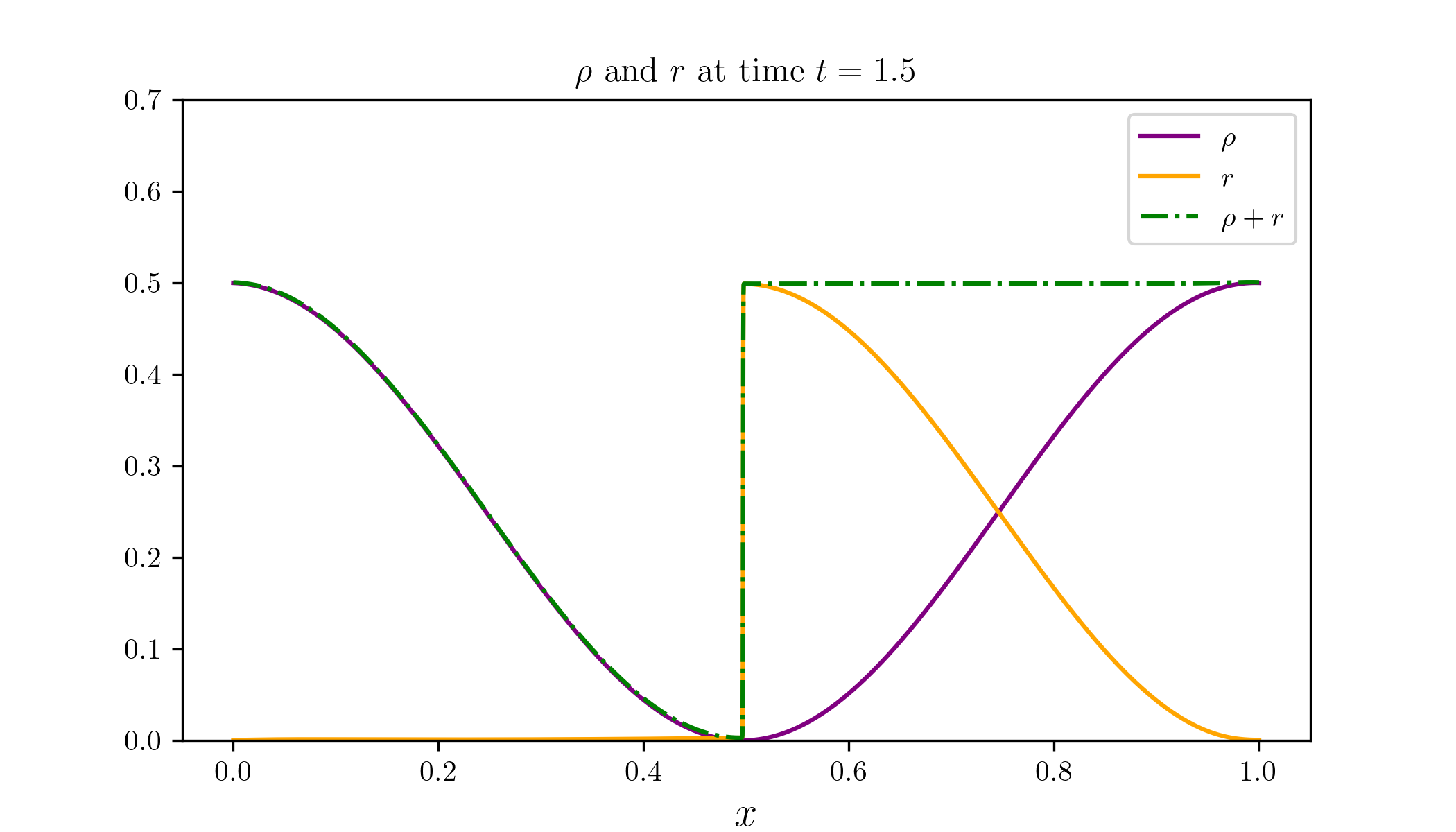}%
  \includegraphics[width=0.333\textwidth,keepaspectratio,trim=25 25 25
  25]{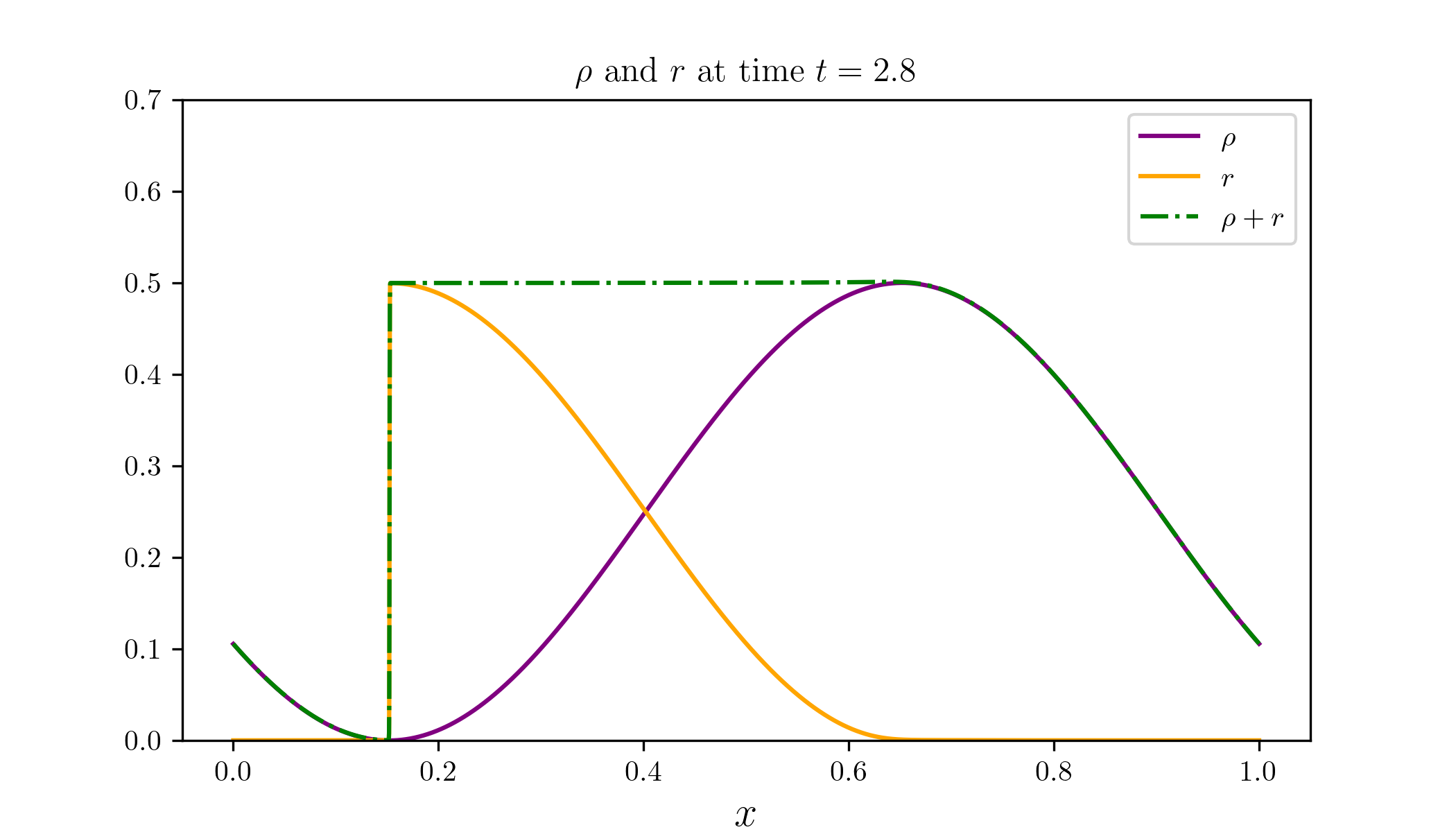}%

  \caption{Solution $(\rho,r)$ to~\eqref{eq:3} with parameters and
    initial datum as in \Cref{rem:r_asymptotic}.  The local component
    $r$, rather than being a traveling wave, converges
    asymptotically-in-time almost uniformly to the wave profile
    $r_\infty$, according to \Cref{thm:Claudia}. Note that in this
    case $\underline{q}_R(x)\equiv 0$ and
    $\overline{q}_R(x)=R-\rho_o(x)$ by~\Cref{rem:ro0} and $r_\infty$
    vanishes on a non trivial interval $I$.  Out of $I$, the sum
    $\rho + r$ is constantly equal to $R$, coherently with
    \Cref{thm:waves} and \Cref{rm:waves}.  The simulation is performed
    using the Lax-Friedrichs method.  }\label{fig:traveling}
\end{figure}

\begin{example}
  \label{rem:r_asymptotic}
  \rm The previous result is consistent with the case of solutions
  $(\rho, r)$ to~\eqref{eq:3} in the sense of \Cref{def:sol} in which
  the nonlocal component $\rho$ is a traveling wave
  whereas the local part $r$ is not, though
  converging asymptotically-in-time to a traveling wave.
  For instance, the choice
  \begin{displaymath}
    \begin{array}{rcl@{\qquad}rcl}
      v_{NL} (r)
       & \coloneq
       & \dfrac12
       & v_L (r)
       & \coloneq
       & \left\{
      \begin{array}{l@{\quad}r@{\,}c@{\,}l}
        1
         & r
         & \in
         & \mathopen]-\infty, 0]
        \\
        1-6r^5+15r^4-10r^3
         & r
         & \in
         & \mathopen]0, 1\mathclose[
        \\
        0
         & r
         & \in
         & [1, +\infty\mathclose[
        \\
      \end{array}
      \right.
      \\
      \rho_o (x)
       & \coloneq
       & \dfrac{1+\sin (2\pi\,x)}{4}
       & r_o (x)
       & \coloneq
       & \dfrac{1-\sin (2\pi\,x)}{8}
    \end{array}
  \end{displaymath}
  satisfies the requirements of Theorem~\ref{thm:WP} and of
  Theorem~\ref{thm:Claudia} with $R=1/2$. See \Cref{fig:traveling}.
\end{example}

\begin{remark}
  \label{rem:ro0}
  \rm In the setting of Theorem~\ref{thm:Claudia}, assume
  $R = \max \Sigma\rho_o$.  Then, for all $x \in \Su$,
  \begin{displaymath}
    \underline{q}_R (x) \equiv 0\,, \quad
    \overline{q}_R (x) \equiv R - \Sigma \rho_o (x)\,,\quad
    f_R = 0\,, \quad
    q_R = 0\,.
  \end{displaymath}
  If moreover $R$ is the unique preimage of $v_L(R)$,
  i.e.~$v_L^{-1}\left(v_L(R)\right) = \{R \}$, then the asymptotic
  profile $r_\infty$ may vanish on a non trivial interval, say $I$,
  and $\Sigma\rho_o$ may well be nonconstant on $I$.  This shows the
  necessity of requirement~\eqref{eq:10} in Theorem~\ref{thm:waves}:
  as soon as $r_\infty$ vanishes on a non trivial interval,
  equality~\eqref{eq:65} may fail, as shown in \Cref{fig:traveling}.
\end{remark}

\section{Proofs}
\label{sec:proofs}

For all $h \in \naturali$, $p \in [1,\infty]$ and
$m\in \naturali \setminus \{0\}$, we identify the spaces
\begin{equation}
  \label{eq:20}
  \small
  \!\!\!
  \begin{array}{@{\!\!\!}r@{\,}c@{\,}l@{\;}r@{\,}c@{\,}l@{\!\!}}
    \L{p} (\Su; \reali^m)
     & \coloneqq
     & \left\{
    f \colon \Su \to \reali^m \colon
    \widehat{f} \in \Lloc{p} (\reali; \reali^m)
    \right\}
     & \norma{f}_{\L{p} (\Su; \reali^m)}
     & \coloneqq
     & \norma{\smash{\widehat{f}}}_{\L{p} ([0,1];\reali^m)}
    \\
    \!
    \W{h}{p}(\Su; \reali^m)
     & \coloneqq
     & \left\{
    f \colon \Su \to \reali^m \colon
    \widehat{f} \in \Wloc{h}{p} (\reali; \reali^m)
    \right\}
     & \norma{f}_{\W{h}{p} (\Su; \reali^m)}
     & \coloneqq
     & \norma{\smash{\widehat{f}}}_{\W{h}{p} ([0,1];\reali^m)}
    \\
    \C{h} (\Su; \reali^m)
     & \coloneqq
     & \left\{
    f \colon \Su \to \reali^m \colon
    \widehat{f} \in \C{h} (\reali; \reali^m)
    \right\}
     & \norma{f}_{\C{h} (\Su; \reali^m)}
     & \coloneqq
     & \norma{\smash{\widehat{f}}}_{\C{h} ([0,1];\reali^m)}
    \\
    \cBV{} (\Su; \reali^m)
     & \coloneqq
     & \left\{
    f \colon \Su \to \reali^m \colon
    \widehat{f} \in \BVloc (\reali; \reali^m)
    \right\}
     & \norma{f}_{\cBV{} (\Su; \reali^m)}
     & \coloneqq
     & \norma{f}_{\L1 (\Su;\reali^m)} {+} \tv (f)
    \\
    \cBV{1} (\Su; \reali^m)
     & \coloneqq
     & \left\{
    f {\in} \W11 (\Su; \reali^m) \colon
    f' \in \cBV{} (\Su; \reali^m)
    \right\}
     & \norma{f}_{\cBV{1}(\Su; \reali^m)}
     & \coloneqq
     & \norma{f}_{\W11(\Su; \reali^m)} {+} \tv(f')
  \end{array}
\end{equation}
If $X,Y$ are any of the spaces above and $I$ is a real interval, we
set
\begin{equation}
  \label{eq:32bis}
  \norma{(f_1,f_2)}_{X\times Y}
  \coloneq
  \norma{f_1}_X + \norma{f_2}_Y
  \quad \mbox{ and } \quad
  \begin{array}{r@{\,}c@{\,}l}
    \norma{f}_{\C0 (I; X)}
     & \coloneq
     & \sup_{t \in I} \norma{f (t)}_X
    \\[6pt]
    \norma{f}_{\L\infty (I; X)}
     & \coloneq
     & \esssup_{t \in I} \norma{f (t)}_X \,.
  \end{array}
\end{equation}
For $f \in \cBV{} (\Su; \reali^m)$, define
\begin{equation}
  \label{eq:44}
  \tv(f)
  \coloneqq
  \tv (\widehat{f}; \mathopen ]0,1]) \,.
\end{equation}
In other words, to compute $\tv(f)$ for
$f \in \cBV{} (\Su; \reali^m)$, first extend $f$ to $\widehat{f}$,
then compute its distributional derivative $\mu$, which is a measure
on $\reali$, pass to its total variation $\norma{\mu}$ and obtain
$\tv(f) = \norma{\mu} (\mathopen]0,1])$.

Given $f \in \L\infty(\Su; \reali), \eta \in \L1(\reali; \reali)$,
define the convolution product $f * \eta$ in $\Su$ as
\begin{equation}
  \label{eq:14}
  f * \eta \coloneqq \widecheck{(\widehat f *
    \eta)}
  \qquad \mbox{ or, equivalently, for } x \in \reali
  \quad
  \widehat{
    (f * \eta)
  }(x)
  =
  {\int_{\reali} \widehat f (x-y) \; \eta (y) \d{y}} \,.
\end{equation}
For completeness, we specify the following definition of use throught:
\begin{equation}
  \label{eq:8}
  \sgn x \coloneq
  \left\{
    \begin{array}{r@{\qquad}r@{\,}c@{\,}l}
      -1
      & x
      & <
      & 0
      \\
      0
      & x
      & =
      & 0
      \\
      1
      & x
      & >
      & 0
    \end{array}
  \right.
\end{equation}
\begin{proofof}{Proposition~\ref{prop:traveling}}
  Observe that by~\eqref{eq:2} and by the regularity of the initial
  datum,
  $(\rho, r) \in \C0\left(\reali_+; \W11(\Su; \reali^k) \times
    \L1(\Su; \reali)\right)$ due to the continuity of translations in
  $\W11(\Su; \reali^k) \times \L1(\Su; \reali)$.

  With reference to the notation in Definition~\ref{def:sol},
  by~\eqref{eq:2} we have
  \begin{eqnarray*}
    v^i (t,x)
    & =
    & v_{NL}^i
      \left(
      \left( \left(\Sigma\widehat{\rho(t)} + \widehat{r (t)}\right) * \eta^i
      \right) (x)
      \right)
      =
      v_{NL}^i
      \left(
      \left( \left(\Sigma\widehat{\rho_o} + \widehat{r_o}\right) * \eta^i
      \right) (x-\lambda\, t)
      \right)
      =
      \lambda \,;
    \\
    w (t,x,u)
    & =
    & v_L\left(\Sigma\widehat{\rho (t)} \, (x) + u \right)
  \end{eqnarray*}
  and note that $w\left( t,x, \widehat{r (t)} (x) \right) = \lambda$.
  Clearly, $(\rho,r)$ is a distributional solution to~\eqref{eq:3}. To
  verify that it satisfies Definition~\ref{def:sol}, compute for all
  $\kappa \in \reali$,
  $\psi \in \Cc1(\mathopen]0, + \infty \mathclose[ \times \reali;
  \reali_+)$
  \begin{eqnarray}
    \nonumber
    &
    & \int_0^{+\infty} \int_{\reali} \modulo{\widehat{r}(t,x) - \kappa} \,
      \partial_t \psi(t,x)
      \d{x}\d{t}
    \\
    \nonumber
    &
    & +
      \int_0^{+\infty} \int_{\reali} \sgn(\widehat{r}(t,x) - \kappa) \,
      \left[
      \widehat{r}(t,x) \, w\left(t,x,\widehat{r} (t,x)\right)
      - \kappa \, w(t,x,\kappa)
      \right] \partial_x \psi(t,x) \d{x}\d{t}
    \\
    \nonumber
    &
    & -
      \int_0^{+\infty} \int_{\reali} \sgn(\widehat{r}(t,x) - \kappa) \,
      \kappa \, \partial_x w(t,x,\kappa) \, \psi(t,x) \d{x}\d{t}
    \\
    \nonumber
    & =
    & \int_0^{+\infty} \int_{\reali} \modulo{\widehat{r}_o(x-\lambda\, t) - \kappa} \,
      \partial_t \psi(t,x)
      \d{x}\d{t}
    \\
    \nonumber
    &
    & + \!\!
      \int_0^{+\infty} \! \int_{\reali}
      \sgn(\widehat{r}_o(x-\lambda\, t) - \kappa)
      \left[
      \lambda \, \widehat{r}_o(x-\lambda\, t) - \kappa
      v_L\left(\Sigma\widehat{\rho_o} (x-\lambda\,t) + \kappa \right)
      \right] \partial_x \psi(t,x) \d{x}\d{t}
    \\
    \nonumber
    &
    & -
      \int_0^{+\infty} \int_{\reali} \sgn(\widehat{r}_o(x-\lambda\, t) - \kappa) \,
      \kappa \,
      v_L'\left( \Sigma\widehat{\rho_o} (x-\lambda\,t)
      + \kappa \right) \,
      \Sigma\widehat{\rho_o}' (x-\lambda\,t) \,
      \psi(t,x) \d{x}\d{t}
    \\
    \label{eq:4}
    & =
    &
      \int_0^{+\infty} \int_{\reali} \modulo{\widehat{r}_o(\xi) - \kappa} \,
      \partial_t \psi(t,\xi+\lambda\, t)
      \d{\xi} \d{t}
    \\
    \label{eq:6}
    &
    &
      +
      \int_0^{+\infty} \int_{\reali} \sgn(\widehat{r}_o(\xi) - \kappa) \,
      \left[
      \lambda \, \widehat{r}_o(\xi) - \kappa \,
      v_L\left( \Sigma\widehat{\rho_o} (\xi) + \kappa \right)
      \right] \partial_x \psi(t,\xi+\lambda\, t) \d{\xi}\d{t}
    \\
    \label{eq:7}
    &
    & -
      \int_0^{+\infty} \int_{\reali} \sgn(\widehat{r}_o(\xi) - \kappa) \,
      \kappa \,
      v_L'\left( \Sigma\widehat{\rho_o} (\xi) + \kappa \right) \,
      \Sigma\widehat{\rho_o}' (\xi) \,\psi(t,\xi+\lambda\, t) \d{\xi}\d{t}
  \end{eqnarray}
  where, in the last equality, we used the change of variable
  $\xi=x-\lambda\, t$.  Introduce the function
  $\Psi(t,\xi) \coloneq \psi(t,\xi+\lambda\, t)$, so that
  \begin{displaymath}
    \partial_t \Psi(t,\xi) =
    \partial_t \psi(t,\xi + \lambda t)
    + \lambda \, \partial_x \psi(t,\xi + \lambda t)
    \quad
    \mbox{and}
    \quad
    \partial_\xi \Psi(t,\xi) = \partial_x \psi(t,\xi+\lambda\, t)
  \end{displaymath}
  and the previous expression becomes
  \begin{eqnarray*}
    \mbox{\eqref{eq:4}}
    & =
    & \int_0^{+\infty} \int_{\reali} \modulo{\widehat{r}_o(\xi) - \kappa} \,
      \left( \partial_t \Psi(t,\xi) - \lambda \, \partial_\xi \Psi (t,\xi) \right)
      \d{\xi} \d{t}
    \\
    \mbox{\eqref{eq:6}}
    & =
    & \int_0^{+\infty} \int_{\reali} \sgn(\widehat{r}_o(\xi) - \kappa) \,
      \left[
      \lambda \, \widehat{r}_o(\xi) - \kappa \,
      v_L\left( \Sigma\widehat{\rho_o} (\xi) + \kappa \right)
      \right] \partial_\xi \Psi(t,\xi) \d{\xi}\d{t}
    \\
    \mbox{\eqref{eq:7}}
    & =
    & -
      \int_0^{+\infty} \int_{\reali} \sgn(\widehat{r}_o(\xi) - \kappa) \,
      \kappa \,
      v_L'\left( \Sigma\widehat{\rho_o} (\xi) + \kappa \right) \,
      \Sigma\widehat{\rho_o}' (\xi) \,\Psi(t,\xi) \d{\xi}\d{t}\,.
  \end{eqnarray*}
  Since
  $\Psi \in \Cc1(\mathopen]0, + \infty\mathclose[ \times \reali;
  \reali_+)$, we get
  $\int_0^{+\infty} \int_{\reali} \modulo{\widehat{r}_o(\xi) - \kappa}
  \, \partial_t \Psi(t,\xi) \d{\xi} \d{t} =0$, so that
  \begin{eqnarray}
    \nonumber
    &
    & \mbox{\eqref{eq:4} $+$~\eqref{eq:6} $+$~\eqref{eq:7}}
    \\
    \nonumber
    & =
    & - \int_0^{+\infty} \int_{\reali} \modulo{\widehat{r}_o(\xi) - \kappa} \,
      \lambda \, \partial_\xi \Psi (t,\xi)
      \d{\xi} \d{t}
    \\
    \nonumber
    &
    & +
      \int_0^{+\infty} \int_{\reali} \sgn(\widehat{r}_o(\xi) - \kappa) \,
      \left[ \lambda \, \widehat{r}_o(\xi) - \lambda \,\kappa \right]
      \partial_\xi \Psi(t,\xi) \d{\xi}\d{t}
    \\
    \nonumber
    &
    & +
      \int_0^{+\infty} \int_{\reali} \sgn(\widehat{r}_o(\xi) - \kappa) \,
      \left[
      \lambda \, \kappa - \kappa \,
      v_L\left( \Sigma\widehat{\rho_o} (\xi) + \kappa \right)
      \right] \partial_\xi \Psi(t,\xi) \d{\xi}\d{t}
    \\
    \nonumber
    &
    & -
      \int_0^{+\infty} \int_{\reali} \sgn(\widehat{r}_o(\xi) - \kappa) \,
      \kappa \,
      v_L'\left( \Sigma\widehat{\rho_o} (\xi) + \kappa \right) \,
      \Sigma\widehat{\rho_o}' (\xi) \,\Psi(t,\xi) \d{\xi}\d{t}
    \\
    \nonumber
    & =
    & \kappa
      \int_0^{+\infty} \int_{\reali} \sgn(\widehat{r}_o(\xi) - \kappa) \,
      \left[
      \lambda  - v_L\left( \Sigma\widehat{\rho_o} (\xi) + \kappa \right)
      \right] \partial_\xi \Psi(t,\xi) \d{\xi}\d{t}
    \\
    \nonumber
    &
    & -
      \kappa
      \int_0^{+\infty} \int_{\reali} \sgn(\widehat{r}_o(\xi) - \kappa) \,
      v_L'\left( \Sigma\widehat{\rho_o} (\xi) + \kappa \right) \,
      \Sigma\widehat{\rho_o}' (\xi) \,\Psi(t,\xi) \d{\xi}\d{t}
    \\
    \label{eq:132}
    & =
    & \kappa
      \int_0^{+\infty} \int_{\reali} \sgn(\widehat{r}_o(\xi) - \kappa) \,
      \partial_\xi \left( \left[
      \lambda  - v_L\left( \Sigma\widehat{\rho_o} (\xi) + \kappa \right)
      \right] \Psi(t,\xi) \right) \d{\xi}\d{t}\,.
  \end{eqnarray}
  For any $t >0$, apply Lemma~\ref{lem:sign} with
  $F (\xi) \coloneq \left[ \lambda - v_L\left( \Sigma\widehat{\rho_o}
      (\xi) + \kappa \right) \right] \Psi(t,\xi)$ and
  $G (\xi) \coloneq \widehat{r}_o(\xi) - \kappa$ to obtain
  \begin{displaymath}
    [\mbox{\eqref{eq:132}}]
    =
    \kappa
    \int_0^{+\infty} \int_{\reali}
    \partial_\xi \left(
      \sgn(\widehat{r}_o(\xi) - \kappa) \,
      \left[
        \lambda  - v_L\left( \Sigma\widehat{\rho_o} (\xi) + \kappa \right)
      \right] \Psi(t,\xi)
    \right)
    \d{\xi}\d{t}
    =
    0
  \end{displaymath}
  since, for all $t \in \reali_+$, the map
  $\left[\xi \mapsto \sgn(\widehat{r}_o(\xi) - \kappa) \, \left( \lambda -
    v_L\left( \Sigma\widehat{\rho_o} (\xi) + \kappa \right) \right)
  \Psi(t,\xi) \right]$ is absolutely continuous with compact support.  The
  proof is completed.
\end{proofof}

\begin{proofof}{\Cref{lem:Mauro}}
  Assume that $c \neq \gamma$. Since $r_o$ is $\C0$ and not constant,
  there exist distinct points $x_m$ and $x_M$ in $\Su$ such that
  $0 \leq r_o (x_m) = \min_{\Su} r_o < \max_{\Su} r_o = r_o
  (x_M)$. Let $x^-_\nu$ and $x^+_\nu$ for $\nu \in \naturali$ be
  sequences such that
  \begin{equation}
    \label{eq:24}
    \!\!\!
    x^-_\nu \neq x^+_\nu, \quad
    r_o (x^-_\nu) = r_o (x^+_\nu) \,,\quad
    \begin{array}{@{}c@{}}
      \lim\limits_{\nu \to +\infty} x^-_\nu = x_M
      \\
      \lim\limits_{\nu \to +\infty} x^+_\nu = x_M
    \end{array}
    \quad \mbox{ so that } \quad
    \begin{array}{@{}l@{\!}}
      \lim\limits_{\nu \to +\infty} r_o (x^-_\nu) = r_o (x_M) \,,
      \\
      \lim\limits_{\nu \to +\infty} r_o (x^+_\nu) = r_o (x_M) \,.
    \end{array}
  \end{equation}
  For any
  $\phi \in \Cc1 (\mathopen]0, + \infty \mathclose[\times \reali;
  \reali)$ such that $\phi (t,x) = \chi (t) \, \psi (x)$, where
  $\psi \in \Cc1 (\reali; \reali)$,
  $\chi \in \Cc1 (\mathopen]0, + \infty \mathclose[; \reali)$ and
  $\int_0^{+\infty}\chi (t)\d{t} = 1$, we have
  \begin{align*}
    \int_0^{+\infty} \!\!\! \int_{\reali}
    \left(
    \widehat{r}(t,x) \,\partial_t \phi(t,x)
    { +}
    \widehat{r}(t,x)  \,
    v_L \left((\widehat{r} {+ }\Sigma\widehat\rho)(t,x)\right) \;
    \partial_x \phi(t,x)
    \right)
    \d{x} \d{t}
    & =0
      \tag*{[\mbox{By Definition~\ref{def:sol}}]}
    \\
    \int_0^{+\infty} \int_{\reali}
    \left(
    \widehat{r}(t,x) \; \partial_t \phi(t,x)
    +
    c \;\widehat{r}(t,x)  \; \partial_x \phi(t,x)
    \right)
    \d{x} \d{t}
    & =0\, .
      \tag*{[\mbox{By~\eqref{eq:22}}]}
  \end{align*}
  Hence, subtracting the above lines, we get
  \begin{align*}
    \int_0^{+\infty} \int_{\reali}
    \widehat{r}(t,x) \left( c -
    v_L \left((\widehat{r} + \Sigma\widehat\rho)(t,x)\right)\right)
    \partial_x \phi(t,x)
    \d{x} \d{t}
    & = 0
    \\
    \!\!\int_{\reali}
    \!\!
    \left[
    \int_0^{+\infty}
    \!
    \widehat{r} (t,x) \left( c {-}  v_L \left(\widehat{r} (t,x) {+ }\Sigma\widehat\rho (t,x)\right)\right)
    \chi (t) \d{t}
    \right]\!
    \partial_x \psi (x)
    \d{x}
    & =
      0
      \tag*{[\mbox{$\phi (t,x) = \chi (t) \, \psi (x)$}]}
    \\
    \int_0^{+\infty}
    \widehat{r} (t,x) \left( c -  v_L \left(\widehat{r} (t,x) + \Sigma\widehat\rho (t,x)\right)\right)
    \chi (t) \d{t}
    & =
      C
      \tag*{[\mbox{For a suitable }$C \in \reali$]}
    \\
    \int_0^{+\infty}
    \left(
    \widehat{r} (t,x)
    \left(
    c -  v_L \left(\widehat{r} (t,x) + \Sigma\widehat\rho (t,x)\right)
    \right)
    -C
    \right)
    \chi (t) \d{t}
    & =
      0
      \tag*{[\mbox{Since }$\textstyle\int\chi=1$]}
    \\
    \widehat{r} (t,x)
    \left(
    c -  v_L \left(\widehat{r} (t,x) + \Sigma\widehat\rho (t,x)\right)
    \right)
    & = C
      \tag*{[\mbox{By Lemma~\ref{lem:variational}}]}
    \\
    \widehat{r}_o (\xi)
    \left(
    c -  v_L \left(\widehat{r}_o (\xi) + \Sigma\widehat{\rho}_o (\xi+ (c-\gamma)t)\right)
    \right)
    & = C,
      \tag*{[\mbox{Set}\,$\xi=x{-}c\, t$, \mbox{use~\eqref{eq:22}}]}
  \end{align*}
  where the last equality holds for all
  $(t,\xi) \in \reali_+ \times \reali$ due to the continuity of
  $\Sigma \rho_o$, $r_o$ and $v_L$.  Choose first $\xi = x^-_\nu$ and
  then $\xi = x^+_\nu$. Then,
  \begin{eqnarray*}
    &
    & r_o (x^-_\nu)
      \left(
      c -  v_L \left(r_o (x^-_\nu)
      + \Sigma\widehat{\rho}_o (x^-_\nu+(c-\gamma)t)\right)
      \right)
    \\
    & =
    & r_o (x^+_\nu)
      \left(
      c -  v_L \left(r_o (x^+_\nu)
      + \Sigma\widehat{\rho}_o (x^+_\nu+(c-\gamma)t)\right)
      \right) \,.
  \end{eqnarray*}
  By~\eqref{eq:24}, $r_o (x^-_\nu) = r_o (x^+_\nu) \neq 0$ for large
  $\nu$, so that for all $t \in \reali_+$
  \begin{displaymath}
    v_L \left(r_o (x^-_\nu)
      + \Sigma\widehat{\rho}_o \left(x^-_\nu+(c-\gamma)t\right)\right)
    =
    v_L \left(r_o (x^+_\nu)
      + \Sigma\widehat{\rho}_o \left(x^+_\nu+(c-\gamma)t\right)\right) \,.
  \end{displaymath}
  By~\eqref{eq:26},\eqref{eq:21} and~\eqref{eq:24},
  \begin{displaymath}
    \mbox{for all } t \in \reali_+
    \qquad
    \Sigma\widehat{\rho}_o \left(x^-_\nu+(c-\gamma)t\right)
    =
    \Sigma\widehat{\rho}_o \left(x^+_\nu+(c-\gamma)t\right)
  \end{displaymath}
  since $c \neq \gamma$ this proves that $\Sigma\widehat{\rho}_o$ is
  $T_\nu$-periodic for all $\nu \in \naturali$, where
  $T_\nu = \modulo{x^+_\nu - x^-_\nu}$. As a consequence,
  $\Sigma\widehat{\rho}_o$ is also
  $\modulo{T_{\nu'} - T_{\nu''}}$-periodic for all $\nu', \nu''$.  The
  fact that $\inf_{\nu', \nu''} \modulo{T_{\nu'} - T_{\nu''}} = 0$
  implies that $\Sigma\widehat{\rho}_o$ is constant, say
  $\Sigma\widehat{\rho}_o (x) = C$ for a suitable $C \in
  \reali$. Then, $\widehat{r}$ is a {Kru\v zkov} solution to the
  scalar (local) equation
  \begin{displaymath}
    \partial_t \widehat{r} + \partial_x \left( \widehat{r} \, v_L (\widehat{r} + C) \right) = 0 \,.
  \end{displaymath}
  Introduce $q_o (x) = \widehat{r}_o (x) - \int_{\Su}r_o (\xi)
  \d\xi$. Clearly,
  $q (t,x) = \widehat{r}_o (x-c\,t) - \int_{\Su}r_o (\xi) \d\xi$
  solves
  \begin{displaymath}
    \partial_t q +
    \partial_x \left(
      \left( q + \int_{\Su}r_o (\xi) \d\xi \right)
      v_L\left( q + \int_{\Su}r_o (\xi) \d\xi  + C  \right)
    \right)
    =0 \,,
  \end{displaymath}
  which we rewrite as $\partial_t q + \partial_x f (q)=0$ with
  $f (q) = \left( q {+ }\int_{\Su}r_o (\xi) \d\xi \right) v_L\left( q
    {+} \int_{\Su}r_o (\xi) \d\xi {+} C \right)$.  Observe
  preliminarily that
  \begin{equation}\label{eq:28}
    \int_{\Su}r_o (\xi)\d{\xi} + C \in \mathopen]0, R\mathclose[\,.
  \end{equation}
  In fact, by contradiction, suppose that
  $\int_{\Su}r_o (\xi)\d{\xi} + C = 0$.  Then, by~\eqref{eq:21} it
  holds
  \begin{align*}
    \int_{\Su}r_o (\xi)\d{\xi}
    +
    \int_{\Su} \Sigma \widehat{\rho}_o (\xi)\d{\xi}
    = 0
    \quad
    \implies
    r_o(\xi) + \Sigma \widehat{\rho}_o (\xi) = 0
    \quad
    \implies
    r_o(\xi) = - C\,,
  \end{align*}
  which is absurd because $r_o$ is not constant.  An analogous
  argument yields $\int_{\Su}r_o (\xi)\d{\xi} + C \neq R$.

  Since
  \begin{displaymath}
    f'' (0)
    =
    2 v_L' \left(\int_{\Su}r_o (\xi)\d{\xi} + C\right)
    +
    \int_{\Su}r_o (\xi)\d{\xi} \;
    v_L'' \left(\int_{\Su}r_o (\xi)\d{\xi} + C\right)
    < 0
  \end{displaymath}
  by~\eqref{eq:26}, \eqref{eq:21} and~\eqref{eq:28}, we can
  apply~\cite[Theorem~1.1]{MR3072229}. This shows that $q$ tends to
  $0$ in $\L{p}(\reali; \reali)$ as $t \to +\infty$, for all
  $p \in [1, +\infty\mathclose[$. Hence, $r$ converges as
  $t \to +\infty$ to a constant function, contradicting the
  assumption.
\end{proofof}

\begin{proofof}{Theorem~\ref{thm:waves}}
  By assumption, for all
  $\phi \in \Cc1(\mathopen]0,+\infty\mathclose[ \times \reali; \reali)$
  \begin{displaymath}
    \int_0^{+\infty} \int_\reali
    \left(
    \widehat{r} \, \partial_t \phi +
    \widehat r \, v_L \left(\widehat r + \Sigma\widehat{\rho}\right) \,
    \partial_x \phi
    \right) \d{x} \d{t} = 0 \,.
  \end{displaymath}
  Choose $\phi \in \Cc1 (\mathopen]0,+\infty\mathclose[ \times \reali; \reali)$ such that
  $\phi(t,x) = \chi(t) \; \psi(x-\lambda \, t)$ with
  $\chi \in \Cc1(\mathopen]0,+\infty\mathclose[; \reali)$, $\int_0^{+\infty}\chi (t) \d{t}=1$ and
  $\psi \in \Cc1(\reali; \reali)$. Then,
  \begin{align*}
    \int_0^{+\infty} \int_\reali
    \big[
     &
      \widehat{r}_o(x- \lambda \, t)
      \left(
      \chi'(t) \, \psi(x-\lambda \, t )
      -
      \lambda \, \chi(t) \, \psi'(x-\lambda \, t )
      \right)
    \\
     & +
      \widehat{r}_o(x- \lambda \, t ) \,
      v_L \left(
      \widehat{r}_o(x-\lambda \, t ) + \Sigma\widehat{\rho}(t,x)
      \right)
      \chi(t) \, \psi'(x-\lambda \, t )
      \big]
    \d{x} \d{t} = 0
  \end{align*}
  and, calling $\xi = x - \lambda \, t$,
  \begin{align*}
    \!\int_0^{+\infty} \!\!\!\int_\reali \!
    \left[
    \widehat{r}_o(\xi)
    \left(
    \chi'(t) \, \psi(\xi )
    {-}
    \lambda \, \chi(t) \, \psi'(\xi)
    \right)\!
    {+}
    \widehat{r}_o(\xi) \, v_L \left( \widehat{r}_o(\xi){+}\Sigma\widehat{\rho}(t,\xi + \lambda \, t) \right)
    \chi(t) \, \psi'(\xi)
    \right]\!\!
    \d{\xi} \d{t}
     & = 0
    \\
    \int_0^{+\infty}
    \int_\reali
    \left[
      \widehat{r}_o(\xi) \, \chi'(t) \, \psi(\xi)
      +
      \widehat{r}_o(\xi) \, \chi(t) \, \psi'(\xi)
      \left(
      -
      \lambda
      +
      v_L \left(\widehat{r}_o(\xi) + \Sigma\widehat{\rho}(t,\xi + \lambda \, t) \right)
      \right)
      \right]
    \d{\xi} \d{t}
     & = 0 .
  \end{align*}
  Since $\chi$ has compact support, the integral of the first summand
  above vanishes, and we get
  \begin{displaymath}
    \int_\reali
    \left[
      \int_0^{+\infty}
      \widehat{r}_o(\xi) \,
      \chi(t)
      \left(
      -
      \lambda
      +
      v_L \left( \widehat{r}_o(\xi) + \Sigma\widehat{\rho}(t,\xi + \lambda \, t) \right)
      \right)
      \d{t}
      \right]
    \psi'(\xi)
    \d{\xi}  = 0,
  \end{displaymath}
  which, due to the arbitrariness of $\psi$, implies that there exists
  $C \in \reali$ such that for a.e.~$\xi \in \reali$,
  \begin{eqnarray*}
    \int_0^{+\infty}
    \widehat{r}_o(\xi)
    \,
    \chi(t)
    \left(
    - \lambda + v_L \left(
    \widehat{r}_o(\xi) + \Sigma\widehat{\rho}(t,\xi + \lambda \, t)
    \right)
    \right)
    \d{t}
    & =
    & C
    \\
    \int_0^{+\infty} \chi(t)
    \left[
      \widehat{r}_o(\xi) \left(
      -\lambda +
      v_L \left(
      \widehat{r}_o(\xi) + \Sigma\widehat{\rho}(t,\xi + \lambda \, t)
      \right)
      \right)
      - C
      \right]
    \d{t}
    & =
    & 0
  \end{eqnarray*}
  and, since $\widehat{r}_o$ is
  right continuous, the above equality
  holds for all $\xi \in \reali$. By an application of
  \Cref{lem:variational}, we obtain that for a.e. $t \in \reali_+$ and
  all $\xi \in \reali$,
  \begin{equation}\label{eq:27}
    \widehat{r}_o(\xi) \left(
    -\lambda +
    v_L \left(
    \widehat{r}_o(\xi) + \Sigma\widehat{\rho}(t,\xi + \lambda \, t)
    \right)
    \right)
    - C = 0 \,.
  \end{equation}
  Since $\Sigma \rho \in \C0(\reali_+; \W11(\Su; \reali))$ and
  $v_L$ is continuous,
  we have that~\eqref{eq:27} holds for
  all $t \in \reali_+$. Use~\eqref{eq:9} at $t = 0$ to obtain $C =
    0$. By~\eqref{eq:10}, since $\widehat{r}_o$ is
  right continuous,
  \begin{displaymath}
    \widehat{r}_o(\xi) + \Sigma \widehat{\rho}(t,\xi)
    \in
    v_L^{-1} \left( \lambda \right)
    \quad \mbox{for all }\, (t,\xi) \in \reali_+\times \reali \,,
  \end{displaymath}
  so that for all $\xi \in \reali$, we have
  \begin{displaymath}
    \Sigma \widehat{\rho}(t,\xi)
    \in
    v_L^{-1} \left( \lambda \right) - \widehat{r}_o(\xi)
    \quad \mbox{for all }\, t \in \reali_+\,.
  \end{displaymath}
  Hence, since $[t \mapsto \Sigma \widehat{\rho}(t,\xi)]$ is
  continuous and by~\eqref{eq:32}, we have that there exists unique
  $s(\xi)\in v_L^{-1} \left( \lambda \right) - \widehat{r}_o(\xi)$
  such that
  $\Sigma \widehat{\rho}(t,\xi) = s(\xi) \quad \mbox{for all }\, t \in
  \reali_+$.  In particular, $\Sigma\widehat\rho (t,\xi)$ is
  independent of $t$, hence $\Sigma \widehat{\rho}$ is a
  $\lambda$-traveling wave.  The proof is completed.
\end{proofof}

\begin{proofof}{Theorem~\ref{thm:Claudia}}
  This proof is inspired by~\cite{Trend-to-steady}. Both tin he
  statement of Theorem~\ref{thm:Claudia} and in the present proof, the
  function $v_L$ is computed exclusively on nonnegative
  values. Therefore, we may assume that for some $R_*>0$, $v_L$ is
  constant for $r \in \mathopen]-\infty , -R_*]$.

  Define $\lambda \coloneq v_L (R)$. By~\eqref{eq:68},
  $\Sigma \widehat{\rho} (t,x) = \Sigma \widehat{\rho}_o (x -
  \lambda\,t)$, hence $\widehat{r}$ solves the Cauchy problem
  \begin{equation}
    \label{eq:25}
    \left\{
    \begin{array}{l@{\qquad}cl}
      \partial_t \widehat r
      +
      \partial_x \left(
      \widehat r \;
      v_L\left( \widehat r + \Sigma \widehat{\rho}_o(x - \lambda t)\right)
      \right) = 0
       & \mbox{in}
       & \reali_+ \times \reali
      \\
      \widehat r(0) =  \widehat{r}_o
       & \mbox{in}
       & \reali \,.
    \end{array}
    \right.
  \end{equation}
  By the above remark on $v_L$, we can apply~\cite[\textbf{(iv)}
    in~Theorem~2.4]{periodic}, proving that the solution $\widehat{r}$
  to problem~\eqref{eq:25} satisfies
  $\widehat{r} (t) \in \BVloc{(\reali; \reali)}$ for all
  $t \in \reali_+$.

  \paragraph{\underline{Step~1:} On $\Su$, problem~\eqref{eq:25} is equivalent to the autonomous problem
    \begin{equation}
      \label{eq:75}
      \left\{
      \begin{array}{@{}ll@{}}
        \partial_t \widecheck{q} + \partial_x f(x, \widecheck{q}) = 0
         & \mbox{in }\reali_+ \times \Su
        \\
        \widecheck{q} (0) = r_o
         & \mbox{in } \reali \,.
      \end{array}
      \right.
      \mbox{ with } \,
      \begin{array}{@{}r@{\,}c@{\,}l@{}}
        q(t, x)
         & \coloneqq
         & \widehat{r}(t, x+\lambda \, t)
        \\
        f(X, Q)
         & \coloneqq
         & Q \left(
        v_L \left( Q + \Sigma\widehat{\rho}_o (X) \right) - \lambda
        \right)
      \end{array}
    \end{equation}
    which admits a solution $\widecheck{q}$ such that
    $\widecheck{q}(t) \in \cBV{}(\Su; \reali)$ for all
    $t \in \reali_+$.}

  Since $\widehat r$ is an entropy solution to~\eqref{eq:25}, then
  $\widecheck{q}$ is an entropy solution to~\eqref{eq:75}, too.
  By~\eqref{eq:26bis} and $\Sigma {\rho}_o \in \C1(\Su; \reali)$,
  \cite[Theorem~2.4]{periodic} yields the uniqueness of entropy
  solutions to~\eqref{eq:75}.  Hence, $\widecheck{q}$ in~\eqref{eq:75}
  is the unique entropy solution to~\eqref{eq:75}.  Moreover,
  $q \in \C0(\reali_+; \Lloc1(\reali; \reali))$ and for all
  $t \in \reali_+$ $q(t)$ is $1$-periodic and in
  $\BVloc(\reali; \reali)$.

  \paragraph{\underline{Step~2:} There exists a unique
    $G \in \C0 ([0,1]; \reali)$ such that%
    \begin{equation} \label{eq:97} \partial_Q f (X,Q) = 0\iff G(X) = Q
      \, \mbox{ and } \, G (X) \in [0, R - \Sigma\rho_o (X)]\
    \end{equation}
    and for all $X \in [0, 1]$
    \begin{equation}
      \label{eq:102}
      \!\!\!
      Q \mapsto f(X, Q) \!\in \C2([0, R-\Sigma \widehat{\rho}_o(X)]; \reali),
      \,
      \begin{array}{c}
        \mbox{is strictly concave, and}
        \\
        \mbox{reaches its maximum at $G\,(X)$.}
      \end{array}
    \end{equation}}
  Introduce the domain
  \begin{equation}
    \label{eq:100}
    \mathcal{D} \coloneq \left\{\left(X, Q\right) \in [0, 1]
    \times \reali_+ \colon Q \le R - \Sigma \widehat{\rho}_o(X)\right\}\,.
  \end{equation}
  and compute
  \begin{eqnarray}
    \label{eq:15}
    \begin{array}{rcl}
      \partial_Q f (X,Q)
       & =
       & v_L \left(Q +\Sigma \widehat{\rho}_o (X)\right)
      - \lambda
      + Q \, v_L'\left(Q + \Sigma \widehat{\rho}_o (X)\right)
      \\
      \partial_X f(X,Q)
       & =
       & Q \,  v_L'\left(Q + \Sigma \widehat{\rho}_o (X)\right) \,
      \Sigma \widehat{\rho}_o'(X)\,.
      \\
      \partial_{QQ}^2 f(X, Q)
       & =
       & 2\, v_L' \left(Q +\Sigma \widehat{\rho}_o (X)\right)
      + Q \, v_L'' \left(Q +\Sigma \widehat{\rho}_o (X)\right) \,.
    \end{array}
  \end{eqnarray}
  By~\eqref{eq:26bis}--\eqref{eq:67} for all
  $(X,Q) \in \pint{\mathcal{D}}$
  \begin{equation}
    \label{eq:96}
    \partial^2_{QQ}f < 0 \,,
  \end{equation}
  $\partial_Q f$ is of class $\C1$ and by~\eqref{eq:26bis}
  \begin{equation}
    \partial_Q f (X,0)
    \geq
    0, \qquad
    \partial_Q f \left(X, R-\Sigma\rho_o (X)\right)
    \leq
    0\,.
  \end{equation}
  Hence, by the Implicit Function Theorem, there exists a unique
  $G \in \C0 ([0, 1]; \reali)$ such that~\eqref{eq:97} holds.  To
  prove the first line in~\eqref{eq:102} use~\eqref{eq:96}, while to
  prove the second line use~\eqref{eq:97}.

  \paragraph{\underline{Step~3:} The function $G$ in~\eqref{eq:97}
    enjoys the properties}
  \begin{enumerate}[label=\textup{\textbf{(G\arabic*)}},
      ref=\textup{\textbf{(G\arabic*)}}]
    \item\label{eq:105}
          $G(X)=0 \iff \Sigma \widehat{\rho}_o (X) = R \iff G(X) = R -
            \Sigma \widehat{\rho}_o (X)$.
    \item\label{eq:101}$\Sigma \widehat{\rho}_o(X)
            \mathrel{\substack{>\\=\\<}} \Sigma \widehat{\rho}_o(X_*) \implies
            G(X) \mathrel{\substack{<\\=\\>}} G(X_*)$.
    \item\label{eq:114}$\Sigma \widehat{\rho}_o(X)
            \mathrel{\substack{>\\=\\<}} \Sigma \widehat{\rho}_o(X_*) \implies
            f(X, G(X)) \mathrel{\substack{<\\=\\>}} f(X_*, G(X_*))$.

    \item \label{item:1} Let $X$ be an isolated stationary point of
          $\Sigma\widehat{\rho}_o$ which is not a point of strict local
          minimum. Then, $\left( X, G (X) \right)$ is a point of saddle for
          $f$.

    \item \label{item:2} Let $X$ be an isolated point of strict local
          minimum of $\Sigma\widehat{\rho}_o$. Then,
          $\left( X, G (X) \right)$ is a point of maximum for $f$.
  \end{enumerate}

  \subparagraph{Proof of~\ref{eq:105}.}  Use the equivalences
  \begin{flalign*}
    G (X) = 0
     & \iff
    \partial_Q f (X,0) = 0
     & [\mbox{By~\eqref{eq:97}}]
    \\
     & \iff
    v_L \left(0 + \Sigma \widehat{\rho}_o (X)\right) = \lambda
     & [\mbox{By~\eqref{eq:15}}]
    \\
     & \iff
    \Sigma \widehat{\rho}_o (X) = R
     & [\mbox{By~\eqref{eq:26bis} and~\eqref{eq:67}}]
    \\
     & \iff
    \left(R - \Sigma \widehat{\rho}_o(X)\right)v_L'(R) = 0
     & [\mbox{By~\eqref{eq:26bis}}]
    \\
     & \iff
    \partial_Q f \left(X, R - \Sigma \widehat{\rho}_o(X)\right) = 0
     & [\mbox{By~\eqref{eq:15}}]
    \\
     & \iff
    G (X) = R - \Sigma \widehat{\rho}_o(X) \,.
     & [\mbox{By~\eqref{eq:97}}]
  \end{flalign*}

  \subparagraph{Proof of~\ref{eq:101}.} First assume that
  $\Sigma \widehat{\rho}_o(X) = \Sigma \widehat{\rho}_o(X_*)$.  Then,
  by~\eqref{eq:75} and~\eqref{eq:97}, one has $G(X) = G(X_*)$.  Assume
  that $\Sigma \widehat{\rho}_o(X)< \Sigma \widehat{\rho}_o(X_*)$.
  Note that by~\eqref{eq:97}
  \begin{equation}\label{eq:113}
    0 \leq G(X_* ) + \Sigma \widehat{\rho}_o(X) < G(X_* ) + \Sigma \widehat{\rho}_o(X_*) \leq R\,.
  \end{equation}
  Since $v_L$ is strictly decreasing and $v_L'$ is non increasing
  on $[0,R]$
  by~\eqref{eq:26bis}, then, using~\eqref{eq:113}, one gets that
  for all $X, X_* \in [0, 1]$
  with $\Sigma \widehat{\rho}_o(X)< \Sigma \widehat{\rho}_o(X_*)$
  \begin{align*}
    0
     & =
    \partial_Q f\left(
    X_*, G(X_*)
    \right)
    \\
     & =
    v_L \left(
    G(X_* ) + \Sigma \widehat{\rho}_o(X_*)
    \right) -\lambda
    + G(X_*) v_L' \left(
    G(X_* ) + \Sigma \widehat{\rho}_o(X_*)
    \right)
    \\
     &
    <
    v_L \left(
    G(X_* ) + \Sigma \widehat{\rho}_o(X)
    \right) -\lambda
    + G(X_*) v_L' \left(
    G(X_* ) + \Sigma \widehat{\rho}_o(X)
    \right)
    \\
     &
    =
    \partial_Q f\left(
    X, G(X_*)
    \right)\,,
  \end{align*}
  which implies that
  \begin{equation}\label{eq:112}
    \partial_Q f\left(
    X, G(X_*)
    \right) > 0.
  \end{equation}
  As the map $Q \mapsto \partial_{q}f(X, Q)$ is strictly decreasing
  by~\eqref{eq:96}, and $\partial_{q}f\left(X, G(X)\right)=0$,
  one concludes from~\eqref{eq:112} that $G(X)> G(X_*)$,
  proving~\ref{eq:101}.

  \subparagraph{Proof of~\ref{eq:114}.}  The case
  $\Sigma \widehat{\rho}_o(X)=\Sigma \widehat{\rho}_o(X_*)$ trivially
  implies~\ref{eq:114} by~\ref{eq:101}.  Suppose that
  $\Sigma \widehat{\rho}_o(X)< \Sigma \widehat{\rho}_o(X_*)$.  First,
  observe that if $\Sigma \widehat{\rho}_o(X_*) = R$, then
  $G(X_*) = 0$ by~\ref{eq:105} and consequently, $f(X_*, G(X_*))=0$.
  Since $v_L$ is strictly decreasing and
  $G(X) \in \mathopen] 0, R-\Sigma \widehat{\rho}_o(X)\mathclose[$
        by~\eqref{eq:97} and~\ref{eq:105}, we obtain
        \begin{align*}
          f \left( X, G(X)
          \right)
          =
          G(X)
          \left(
          v_L \left(
          G(X) + \Sigma \widehat{\rho}_o(X)
          \right) - \lambda
          \right)
          > 0 = f(X_*, G(X_*))\,.
        \end{align*}
        In the complementary case
      $\Sigma \widehat{\rho}_o(X_*) < R$,
        by~\ref{eq:105} we have $G(X_*) >0$, so that
        \begin{align*}
          f \left( X, G(X)
          \right)
           &
          =
          G(X)
          \left(
          v_L \left(
          G(X) + \Sigma \widehat{\rho}_o(X)
          \right) - \lambda
          \right)
          \\
           &
          \geq
          f(X, G(X_*))
          \tag*{[\mbox{By~\eqref{eq:102}}]}
          \\
           &
          =
          G(X_*)
          \left(
          v_L \left(
          G(X_*) + \Sigma \widehat{\rho}_o(X)
          \right) - \lambda
          \right)
          \\
           &
          >
          G(X_*)
          \left(
          v_L \left(
          G(X_*) + \Sigma \widehat{\rho}_o(X_*)
          \right) - \lambda
          \right)
          \tag*{[\mbox{Since $v_L$ is strictly decreasing}]}
          \\
           &
          =
          f \left( X_*, G(X_*)
          \right)\, ,
        \end{align*}
        proving~\ref{eq:114}.

        \subparagraph{Proof of~\ref{item:1}.} Assume first that $X$ is a
        strict local maximum of $\Sigma \widehat{\rho}_o$, there exists
      $\varepsilon > 0$ such that
      $\Sigma \widehat{\rho}_o(X_*) < \Sigma \widehat{\rho}_o(X)$ for all
      $X_* \in \mathopen] X - \varepsilon, X + \varepsilon \mathclose[
      \setminus \{X\}$.  An application of~\ref{eq:114} yields that
      $f\left(X_*, G(X_*)\right) > f\left(X, G(X)\right)$, proving that
      $X$ is strict local minimum of the map
      $[X \mapsto f \left(X, G(X)\right)]$.  On the other hand, recall
        by~\eqref{eq:102} that the map $[Q \mapsto f(X, Q)]$ is strictly
        concave and reaches its maximum at $G(X)$.  Since we have found a
        direction along which $f$ is strictly increasing and another one
        along which $f$ is strictly decreasing at $(X, G(X))$, we conclude
        that $(X, G(X))$ is a saddle of $f$.

        If $X$ is not a point of local maximum, an analogous reasoning applies.

        \subparagraph{Proof of~\ref{item:2}.} The point $X$ is a strict
        local minimum of the map $\Sigma \widehat{\rho}_o$, hence there
        exists $\varepsilon > 0$ such that
      $\Sigma \widehat{\rho}_o(X_*) > \Sigma \widehat{\rho}_o(X)$ for all
      $X_* \in \mathopen] X - \varepsilon, X + \varepsilon \mathclose[
    \setminus \{X\}$.  Use~\ref{eq:114} to get
  $f\left(X_*, G(X_*)\right) < f\left(X, G(X)\right)$, proving that
  $X$ is strict local maximum of the map
  $[X \mapsto f \left(X, G(X)\right)]$.  Fix any
  $(X_*, Q_*) \in B_\varepsilon(X, G(X)) \setminus \{(X, G(X))\}$,
  where $B_\varepsilon(X, G(X))$ is the open ball of radius
  $\varepsilon$ centered at $(X, G(X))$.  If $X_*=X$ and
  $Q_* \neq G(X)$, then $f(X_*, Q_*) < f(X, G(X))$ by~\eqref{eq:102}.
  Otherwise, we have that
  \begin{align*}
    \!\!
    f(X_*, Q_*)
     &
    \leq
    f(X_*, G(X_*))
    \tag*{
      [\mbox{$G(X_*)$ is strict local maximum of $[Q \mapsto f(X_*, Q)]$
            by~\eqref{eq:102}}]}
    \\
     &
    <
    f(X, G(X)),
    \tag*{[\mbox{$X$ is strict maximum of $[X \mapsto f(X, G(X))]$}]}
  \end{align*}
  proving that $(X, G(X))$ is a strict local maximum of $f$.

  \paragraph{\underline{Step~4:} There exist functions
    $\underline{q}_R, \overline{q}_R$ uniquely defined
    by~\ref{item:43} in theorem~\ref{thm:Claudia}. Moreover, $\underline{q}_R, \overline{q}_R$ are continuous and satisfy
    \begin{gather}
      \label{eq:158}
      f\left( X, \underline{q}_R (X) \right) = f_R
      \qquad
      \mbox{and} \qquad
      f\left( X, \overline{q}_R (X) \right) = f_R
      \\
      \label{eq:107}
      Q \in \mathopen]\underline{q}_R (X),\overline{q}_R (X)\mathclose[
      \iff
      f (X,Q) > f_R
      \\
      \label{eq:106}
      \underline{q}_R(X) =\overline{q}_R(X)
      \iff
      f\left(X, G(X)\right) = f_R
      \iff
      \Sigma  \widehat{\rho}_o(X) = \max_{\Su} \Sigma {\rho}_o\,.
    \end{gather}}
  Fix $X_o \in \mbox{argmax}_{[0,1]} \Sigma \widehat{\rho}_o$,
  $q_R \coloneqq G(X_o)$ and set $f_R \coloneq f(X_o, q_R)$.
  By~\ref{eq:101} and~\ref{eq:114}, the values $q_R$ and $f_R$ are
  independent of the choice of
  $X_o \in \mbox{argmax}_{[0,1]} \Sigma \widehat{\rho}_o$.

  The definition of $f$ in~\eqref{eq:75} yields
  \begin{align*}
    f_R = 0
     &
    \iff q_R = 0
    \quad
    \mbox{or}
    \quad
    v_L\left(q_R + \Sigma \widehat{\rho}_o(X_o)\right) = \lambda \,
    \tag*{[By~\eqref{eq:75}]}
    \\
     &
    \iff G(X_o) = 0
    \quad
    \mbox{or}
    \quad
    G(X_o) + \Sigma \widehat{\rho}_o(X_o) = R \,
    \tag*{[By~\eqref{eq:26bis}]}
    \\
     &
    \iff
    \Sigma \widehat{\rho}_o(X_o) = R \,,
    \tag*{[By~\ref{eq:105}]}
  \end{align*}
  so that, using~\eqref{eq:67} and~\eqref{eq:26bis}, one gets
  \begin{equation}\label{eq:103}
    f_R \geq 0
    \qquad
    \mbox{ and }
    \qquad
    f_R = 0 \iff
    \max_{\Su} \Sigma {\rho}_o = R\, .
  \end{equation}
  An application of~\ref{eq:114} with $X_*=X_o$ yields
  that for all $X \in [0,1]$,
  \begin{align}\label{eq:104}
    f\left(X, G(X)\right) \geq f_R
    \qquad
    \mbox{ and }
    \qquad
    f\left(X, G(X)\right) = f_R
     &
    \iff
    \Sigma  \widehat{\rho}_o(X) = \max_{\Su } \Sigma{\rho}_o\,.
  \end{align}
  \noindent
  By~\eqref{eq:102} and~\eqref{eq:104}, there exist unique
  $\underline{q}_R, \overline{q}_R \colon [0,1] \to \reali$
  such that for all $X \in [0,1]$,
  \begin{gather}
    \label{eq:115}
    0
    \leq \underline{q}_R(X)
    \leq
    G(X)
    \leq \overline{q}_R(X)
    \leq
    R - \Sigma \widehat{\rho}_o(X)
    \\
    \nonumber
    f\left( X, \underline{q}_R (X) \right) = f_R
    \qquad
    \mbox{and} \qquad
    f\left( X, \overline{q}_R (X) \right) = f_R \,.
  \end{gather}
  Moreover, for all $(X,Q) \in \mathcal{D}$, \eqref{eq:102},
  \eqref{eq:104}, \eqref{eq:115} and~\eqref{eq:158} imply
  that~\eqref{eq:107} and~\eqref{eq:106} hold.
  By~\eqref{eq:102}, we obtain that
  $\underline{q}_R, \overline{q}_R \in \C0([0,1]; \reali)$,
  completing the proof of \textbf{\underline{Step~4}}.

  \bigskip

  The solution to~\eqref{eq:75} is deeply related to the
  characteristics~\cite[Chapter~11.9]{dafermos}, i.e., to the solutions
  to the Hamiltonian system of ordinary differential equations
  \begin{equation}
    \label{eq:70}
    \!\!\!\!\!
    \left\{
    \begin{array}{@{\,}l@{}}
      \dot X = \partial_Q f (X,Q)
      \\
      \dot Q  = - \partial_X f(X,Q).
    \end{array}
    \right.
  \end{equation}

  \paragraph{\underline{Step~5:} The curves $Q = \underline{q}_R (X)$
    and $Q = \overline{q}_R (X)$ contain heteroclinic orbits
    of~\eqref{eq:70}, the region
    $\mathcal{B} \coloneqq \left\{ (X,Q) \in \mathcal{D} \colon Q \in [\underline{q}_R (X),
        \overline{q}_R (X)] \right\}$ is invariant for~\eqref{eq:70}
    and contains a finite number of homoclinic/heteroclinic orbits
    of~\eqref{eq:70}.
    Moreover,
    the only possible orbits in the region
    $\mathcal{B}$ are
    either periodic or lie in
    heteroclinic/homoclinic curves.
  }

  Since $f$ is constant along each solution to~\eqref{eq:70},
  using~\eqref{eq:158} and~\eqref{eq:107},
  it holds that the set $\mathcal{B}$ is invariant for~\eqref{eq:70}.

  We tackle the analysis of the equilibria
  of~\eqref{eq:70}, which are strictly related to the homoclinic and
  heteroclinic orbits of~\eqref{eq:70}.
  By~\eqref{eq:26bis},~\eqref{eq:75},~\eqref{eq:97}
  and~\eqref{eq:70}, in $\mathcal{D}$
  \begin{align*}
    \dot Q = 0
    \iff Q = 0 \mbox{ or } \Sigma \widehat{\rho}_o' (X) = 0
    \quad \mbox{and}
    \quad
    \dot X = 0
    \iff Q= G(X) \,.
  \end{align*}
  Hence, the set of the points of equilibrium
  of~\eqref{eq:70} is
  \begin{displaymath}
    \mathcal{E}
    \coloneqq
    \left\{
    \left(X, G(X)\right) \in \mathcal{D}
    \colon X \in \mathcal{S}
    \right\}
    \quad
    \mbox{where}
    \quad
    \mathcal{S}
    \coloneqq
    \left\{
    X \in [0,1]
    \colon \Sigma \widehat{\rho}_o'(X) = 0
    \right\}.
  \end{displaymath}
  By~\eqref{eq:5}, the set $\mathcal{S}$ is a finite union of closed
  intervals.
  By~\eqref{eq:158}, the curves $Q = \underline{q}_R (X)$
  and $Q = \overline{q}_R (X)$ are invariant for~\eqref{eq:70}.
  Since $(0, G(0))$ and $(1, G(1))$ are equilibrium points,
  and
  \begin{displaymath}
    \begin{array}{lll}
      \dot X > 0
       &
      \iff
       &
      \{ (X, Q) \in \mathcal{D}
      \colon \,G(X) > Q
      \}
      \\
      \dot Q > 0
       &
      \iff
       &
      \{ (X, Q) \in \mathcal{D}
      \colon \, \Sigma \widehat{\rho}_o'(X) > 0\, ,
      Q > 0
      \}
    \end{array}
  \end{displaymath}
  we obtain that $Q = \underline{q}_R (X)$ and
  $Q = \overline{q}_R (X)$ are heteroclinic curves for~\eqref{eq:70},
  possibly connecting additional points of equilibrium $(X, G(X))$
  with $ X \in \rm{argmax}_{[0,1]}\Sigma \widehat{\rho}_o$
  by~\eqref{eq:104}.

  Fix $\bar X \in \mathcal{S}$ and let $\mathcal{L}$ be either
  $\mathcal{L} \coloneqq \mathopen]\bar X - \varepsilon, \bar
  X\mathclose[$ or
  $\mathcal{L} \coloneqq \mathopen]\bar X, \bar X +
  \varepsilon\mathclose[$, with $\varepsilon > 0$.  By~\eqref{eq:5},
  choosing $\varepsilon$ sufficiently small, exactly one of the
  following holds on $\mathcal{L}$:
  \begin{enumerate}[label=\textup{\textbf{($\mathcal{L}$\arabic*)}},
    ref=\textup{\textbf{($\mathcal{L}$\arabic*)}}]
  \item\label{item:4}$\Sigma \widehat{\rho}_o (X) < \Sigma
    \widehat{\rho}_o (\bar X)$ for all $X \in \mathcal{L}$.
    By~\ref{eq:114}, this implies that
    $f(X, G(X)) > f(\bar X, G(\bar X))$ for $X \in \mathcal{L}$, and
    by~\eqref{eq:102} there exist exactly two functions
    $\underline q, \overline q \colon \mathcal{L} \to \reali$, such
    that $\underline q (X) < \overline q(X)$ and
    $f(X, \underline q(X)) = f(X, \overline q(X)) = f(\bar X, G(\bar
    X))$ for all $X \in \mathcal{L}$.  Hence, on
    $\left( \mathcal{L} \times \reali_+ \right) \cap \mathcal{D}$, the
    level set of $f$ through the point $(\bar X, G(\bar X))$ is
    composed by two arcs only meeting at $(\bar X, G(\bar X))$.

  \item\label{item:5}
    $\Sigma \widehat{\rho}_o (X) >\Sigma \widehat{\rho}_o (\bar X)$
    for all $X \in \mathcal{L}$.  Since by~\ref{eq:114}
    $f(X, G(X)) < f(\bar X, G(\bar X))$ for all $X \in\mathcal{L}$,
    then, by~\eqref{eq:102}, the level set of $f$ passing through
    $(\bar X, G(\bar X))$ does not intersect the set
    $\left( \mathcal{L} \times \reali_+ \right) \cap \mathcal{D}$.
  \item\label{item:6}
    $\Sigma \widehat{\rho}_o (X) = \Sigma \widehat{\rho}_o (\bar X)$
    for all $X \in \mathcal{L}$.  Using~\ref{eq:114} the level set of
    $f$ through the point $(\bar X, G(\bar X))$ on
    $\left( \mathcal{L} \times \reali_+ \right) \cap \mathcal{D}$ is
    composed by a single arc.
  \end{enumerate}
  The homoclinic/heteroclinic curves in $\mathcal{D}$ are connected to
  the equilibrium points of~\eqref{eq:70}.  We prove that the number
  of such curves is bounded by $4 \, \sharp(\partial \mathcal{S})$,
  $\sharp (E)$ being the cardinality of $E$. Indeed:
  \begin{itemize}
  \item If $\bar X \in \mathcal{S}$ is isolated and strict maximum of
    $\Sigma \widehat{\rho}_o$, then~\ref{item:4} holds on both sides
    of $\bar X$ and, consequently, the level set of $f$ passing
    through $(\bar X, G(\bar X))$ is locally composed by four arcs
    meeting at $(\bar X, G(\bar X))$.  By~\ref{item:1},
    $(\bar X, G(\bar X))$ is a saddle point both for $f$ and for
    system~\eqref{eq:70} by~\cite[Theorem~2,
    Section~2.14]{zbMATH01581767}.  In that case, at most $4$
    homoclinic/heteroclinic curves are related to
    $(\bar X, G(\bar X))$.
  \item If $\bar X\in \mathcal{S}$ is isolated and strict minimum of
    $\Sigma \widehat{\rho}_o$, then~\ref{item:5} holds on both sides
    of $\bar X$.  By~\ref{item:2}, $(\bar X, G(\bar X))$ is a point of
    maximum for $f$ and a center for system~\eqref{eq:70}
    by~\cite[Theorem~2, Section~2.14]{zbMATH01581767}.  In that case,
    no homoclinic/heteroclinic curve is related to
    $(\bar X, G(\bar X))$.

  \item If $\bar X \in \mathcal{S}$ is isolated, neither maximum nor
    minimum for $\Sigma \widehat{\rho}_o$, then~\ref{item:4} holds on
    one side of $\bar X$ and~\ref{item:5} on the other.  Hence, the
    level set of $f$ passing through $(\bar X, G(\bar X))$ is locally
    composed by two arcs meeting at $(\bar X, G(\bar X))$.  In that
    case, at most $2$ homoclinic/heteroclinic curves are related to
    $(\bar X, G(\bar X))$.

  \item If $\bar X \in \mathcal{S}$ belongs to the interior of a
    maximal interval $[X_1, X_2]$ of $\mathcal{S}$, since the vector
    field of the system~\eqref{eq:70} is purely horizontal in the
    strip $[X_1,X_2] \times \reali_+$, no non trivial trajectory
    of~\eqref{eq:70} (and consequently no homoclinic nor heteroclinic
    curve) can approach the equilibrium point $(\bar X,G(\bar X))$.
  \item If $\bar X \in \mathcal{S}$ is the endpoint of a maximal
    interval $[X_1,X_2]$ of $\mathcal S$ (say $\bar X = X_1$),
    \ref{item:6} holds on the interior side; on the outer side either
    case~\ref{item:4} or case~\ref{item:5} holds.  Either way, the
    level set of $f$ passing through $(\bar X, G(\bar X))$ is locally
    composed by finitely many arcs and at most $2$
    homoclinic/heteroclinic curves are related to
    $(\bar X, G(\bar X))$.
  \end{itemize}
  The above analysis proves that $\mathcal{B}$ contains finitely many
  homoclinic/heteroclinic orbits of system~\eqref{eq:70}.  The
  homoclinic/heteroclinic curves are separatrices confining the
  periodic orbits, by the standard classification of the phase
  portrait of planar Hamiltonian systems.

  \begin{figure}
    \centering \centering
    \includegraphics[width=\textwidth]{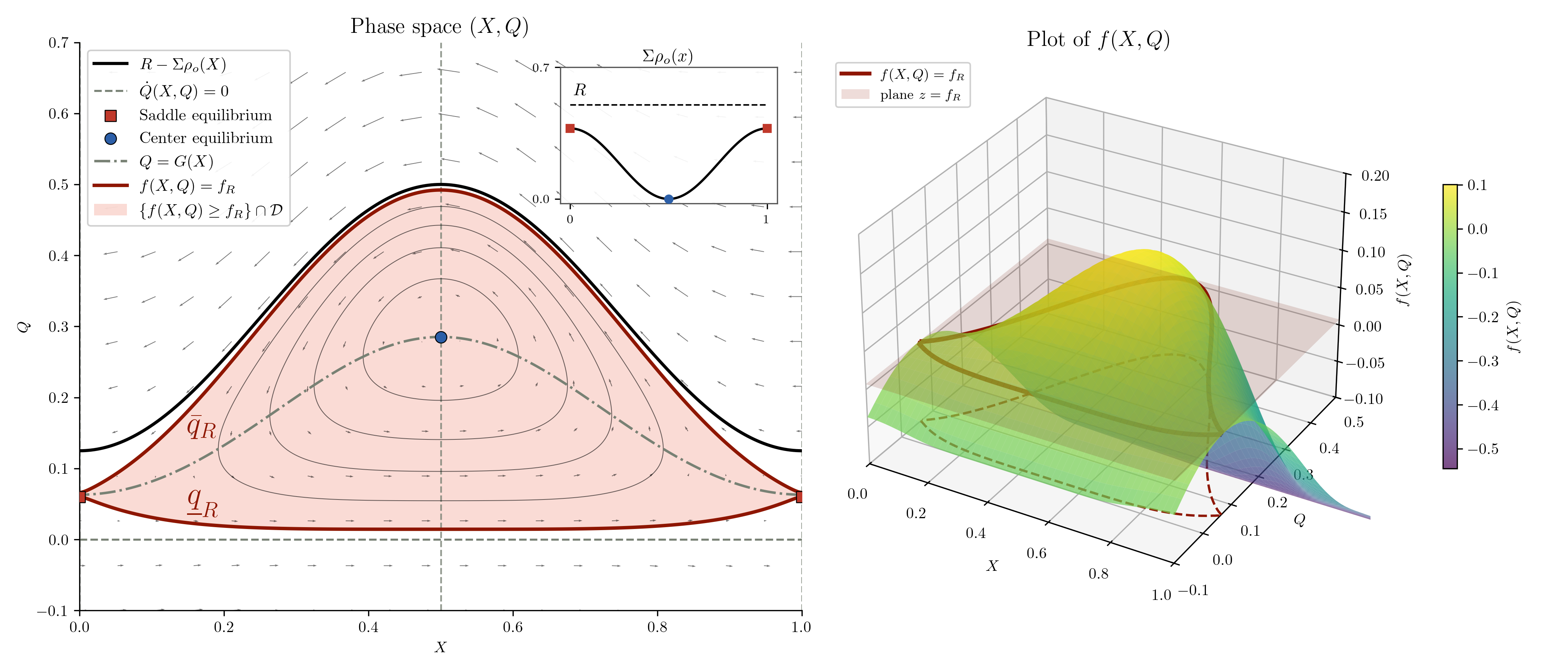}
    \caption{\textbf{\underline{Step~5}} with
      $\Sigma\rho_o$ having a unique point of global maximum. Left,
      the phase plane $(X,Q)$ and, right, the space $(X,Q,f)$ related
      to~\eqref{eq:70}. The equilibria are one center and two saddle
      points. All orbits confined between $\underline{q}_R$ and
      $\overline{q}_R$ are periodic.}
    \label{fig:1a}
  \end{figure}

  \paragraph{\underline{Step~6:}
    It holds that
    \begin{equation}\label{eq:87}
      \lim_{t \to + \infty }
      \,\sup_{x \in [0,1]}
      \,\min_{\bar f \in \mathcal{F}}
      \modulo{
        f\left(x, q(t,x)\right) - \bar f
      }
      =0\, ,
    \end{equation}
    where $\mathcal{F} \coloneqq \left\{
      f\left(X, G(X) \right) \colon X \in \partial \mathcal{S},
      X \textrm{ not a strict minimum for } \Sigma \widehat{\rho}_o
      \right\}$.
  }

  Observe that $\mathcal{F}$ is a finite set, since
  $\partial \mathcal{S}$ is finite.
  The functions $\underline{q}_R$ and $\overline{q}_R$
  in \textbf{\underline{Step~4}} are stationary solutions
  to the equation in \eqref{eq:75} and along them the
  flux $f$ is constantly equal to $f_R$.
  Hence, by~\cite[Theorem~3]{Kruzkov1970},
  since
  $\underline{q}_R \leq \widehat{r}_o \leq \overline{q}_R$
  by~\eqref{eq:138},
  then for a.e. $(t,x) \in \reali_+ \times [0,1]$,
  $
    \underline{q}_R \leq q(t, x) \leq \overline{q}_R
  $.
  Using the fact that,
  for all $t\in \reali_+$,
  $q(t) \in \BV([0,1]; \reali)$ and
  right continuous, we have
  \begin{equation}\label{eq:119}
    \underline{q}_R \leq q(t, x) \leq \overline{q}_R
    \quad
    \mbox{ for all }
    \,
    (t, x) \in \reali_+ \times [0,1]\,.
  \end{equation}
  Fix $(t, x) \in \reali_+ \times [0,1]$ and consider the maximal
  backward generalized characteristics propagating from $(t,x)$, which
  we denote as $s \mapsto X(s; t, x)=$, for $s \in [0,t]$.  Moreover,
  call $Q(s; t, x) \coloneqq q\left(s, X(s; t, x)\right)$, for
  $s \in [0,t]]$.  The function
  $\left(X(\cdot; t,x), Q(\cdot; t,x)\right)$ is the solution to the
  system~\eqref{eq:70} on $[0,t]$ with final data
  $ \left( X(t; t,x), Q(t; t,x) \right) = ( x, q(t,x) ) $,
  see~\cite[Theorem~11.9.3]{dafermos}. Introduce the support of the
  corresponding orbit
  \begin{displaymath}
    \mathcal{O}\left(x,q(t,x)\right) \coloneqq
    \left\{
    \left(
    X(s; t,x), Q(s; t,x)
    \right)
    \colon
    s \in [0, t]
    \right\}
  \end{displaymath}
  and
  $\mathcal{O}\left(x,q(t,x)\right) \subset f^{-1}\left(f\left(x,
      q(t,x)\right)\right)$ since $f$ is constant along the solutions
  to~\eqref{eq:70}.

  Suppose that $f\left(x, q(t,x)\right) \notin \mathcal{F}$.
  Assume that $(x, q(t,x))$ coincides with
  a center of~\eqref{eq:70}.
  Hence, the associated maximal backward generalized characteristics $X(s, t,x) = x$  for all $s \in [0,t]$.
  Since $q(t)$ is right continuous
  and $f$ is continuous,
  there exists $\varepsilon >0$ such that the value
  $f\left(x+\varepsilon, q(t,x+\varepsilon)\right) \notin \mathcal{F}$
  and $(x+\varepsilon, q(t, x+\varepsilon ))$ does not coincide with a center.
  Hence, the orbit $\mathcal{O}(x+\varepsilon, q(t, x+\varepsilon))$
  is periodic with period $T_\varepsilon$.
  If $T_\varepsilon > t$, the characteristics $X(\cdot, t,x)$ and
  $X(\cdot, t,x+\varepsilon)$ collide at a certain
  $s_* \in [0,t]$, which is absurd by~\cite[Theorem~11.9.4]{dafermos}.

  Therefore, we can assume that
  $(x, q(t,x))$ does not
  coincide with a center of~\eqref{eq:70}.
  By \textbf{\underline{Step~5}}, $\mathcal{O}\left(x, q(t,x)\right)$
  is contained in a periodic orbit
  confined by a heteroclinic/homoclinic curves.
  Its period $T$ can be computed
  as
  \begin{displaymath}
    T \coloneqq \oint_{\{f= c\}} \frac{d\sigma}{|\nabla f|}
    \quad \mbox{with $c = f\left(x, q(t,x)\right)$}.
  \end{displaymath}
  Since $\nabla f$ is continuous, $\nabla f \neq 0$ along
  $\{f=c\}$
  and the curve
  $\{f=c\}$ continuously depends on $c$,
  we have that also
  $T$ continuously depends on $c$.
  Moreover, the more $f\left(x, q(t,x)\right)$ is close to some $\bar f$ in $\mathcal{F}$,
  the larger $T$ is, and $T\to + \infty$ as
  $f\left(x, q(t,x)\right) \to \bar f$.

  \medskip
  We suppose that $t > T$ and claim that this leads to a
  contradiction.

  Since $q(t)$ is right continuous and $f$ is continuous, there exists
  $\varepsilon >0$ such that the value
  $f\left(x+\varepsilon, q(t,x+\varepsilon)\right)$ is arbitrarily
  close to $f\left(x, q(t,x)\right)$ and does not belong to
  $\mathcal{F}$.  Consequently,
  $\mathcal{O}\left(x+\varepsilon,q(t,x+\varepsilon)\right)$ is a
  periodic orbit as well, arbitrarily close to
  $\mathcal{O}\left(x, q(t,x)\right)$.  Also, its period, call it
  $T_\varepsilon$, is arbitrarily close to $T$.  Hence, we may assume
  $t > \max\{T, T_\varepsilon\}$.

  In the case
  $f\left(x+\varepsilon, q(t,x+\varepsilon)\right) > f\left(x,
    q(t,x)\right)$ (respectively if
  $f\left(x+\varepsilon, q(t,x+\varepsilon)\right) < f\left(x,
    q(t,x)\right)$), then the orbit
  $\mathcal{O}\left(x+\varepsilon,q(t,x+\varepsilon)\right)$ is
  confined by $\mathcal{O}\left(x, q(t,x)\right)$ (respectively, it
  confines $\mathcal{O}\left(x,q(t,x)\right)$).  Then, if
  $f\left(x+\varepsilon, q(t,x+\varepsilon)\right) = f\left(x,
    q(t,x)\right)$, the two orbits overlap.  See \Cref{fig:4}, left
  images, for a visual illustration.

  \begin{figure}[h!]
    \centering
    \includegraphics[width=\textwidth]{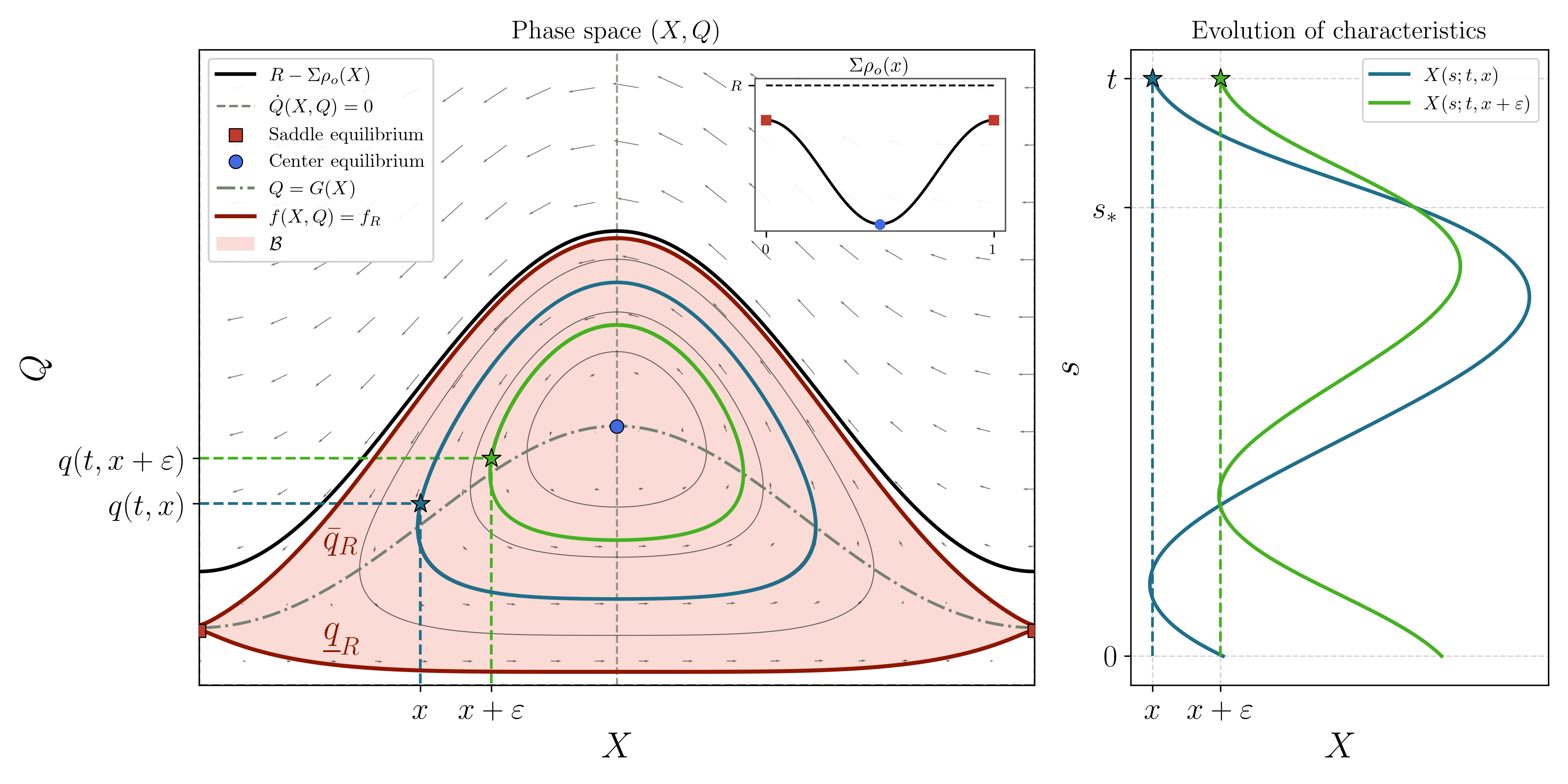}%
    \\
    \includegraphics[width=\textwidth]{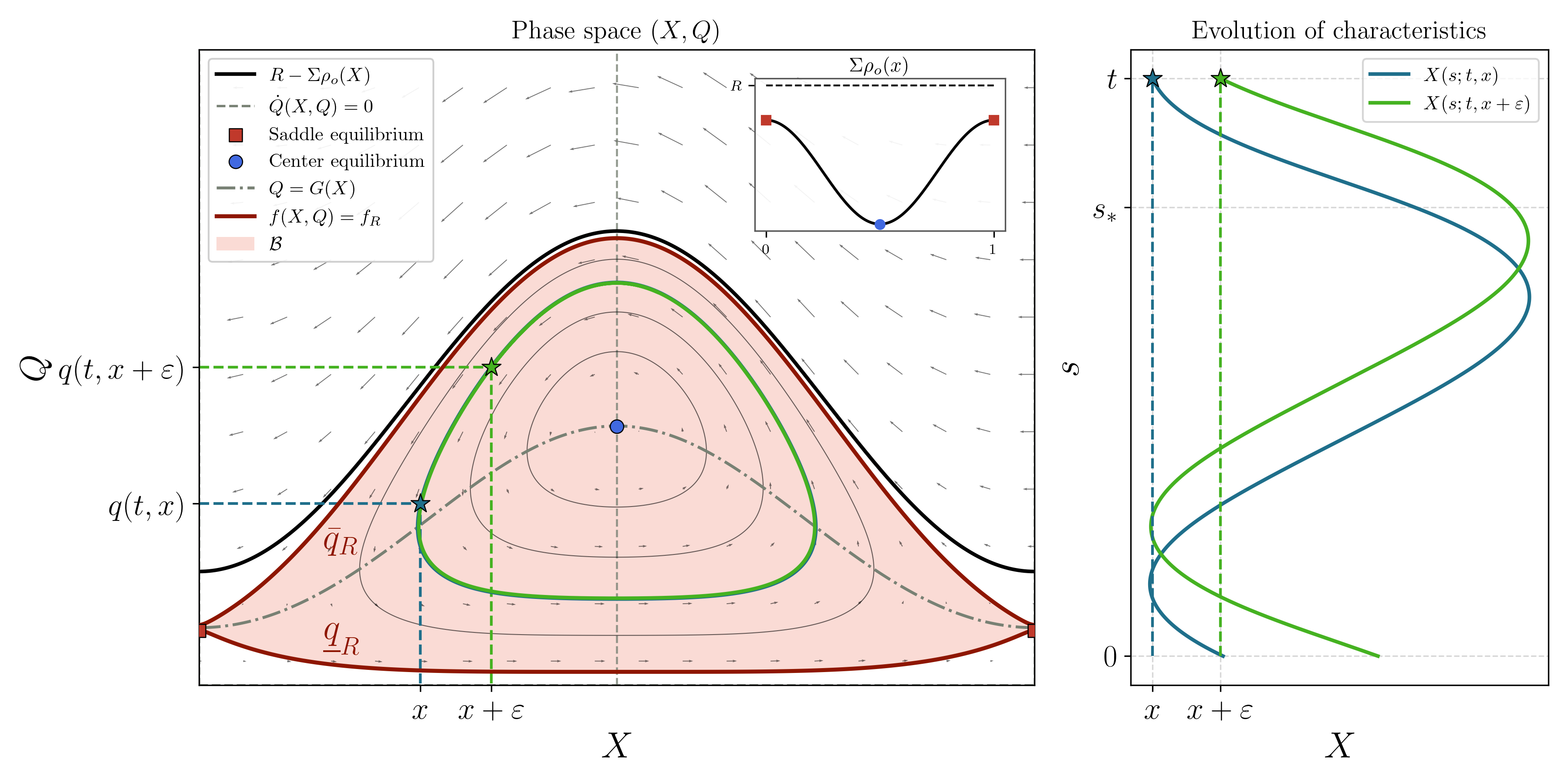}%
    \caption{Maximal backward generalized characteristics, above when
      $f\left(x+\varepsilon, q(t,x+\varepsilon)\right) > f\left(x,
        q(t,x)\right)$ and, below, when
      $f\left(x+\varepsilon, q(t,x+\varepsilon)\right) = f\left(x,
        q(t,x)\right)$.  On the left, in the $(X,Q)$ plane, the
      periodic orbits $\mathcal{O}\left(x, q(t,x)\right)$ (blue) and
      $\mathcal{O}\left(x+\varepsilon, q(t,x+\varepsilon)\right)$
      (green).  On the right, in the $(X,s)$ plane the characteristics
      $s \mapsto X(s; t,x)$ (blue) and
      $s \mapsto X(s; t,x+\varepsilon)$ (green) for $s \in [0,t]$.  In
      both cases, the maximal backward generalized characteristics
      intersect at $s_*$.}
    \label{fig:4}
  \end{figure}

  Suppose without loss of generality that
  $f\left(x+\varepsilon, q(t,x+\varepsilon)\right) \geq f\left(x,
    q(t,x)\right)$.  The maximal backward generalized characteristics
  $X(\cdot; t, x)$ and $X(\cdot; t, x+\varepsilon)$ are continuous
  periodic functions such that
  \begin{displaymath}
    \min_{s \in [0,t]} X(s; t,x)
    \leq
    \min_{s \in [0,t]} X(s; t,x+\varepsilon)
    \leq
    \max_{s \in [0,t]} X(s; t,x+\varepsilon)
    \leq
    \max_{s \in [0,t]} X(s; t,x).
  \end{displaymath}
  Since $t > \max\{T, T_\varepsilon\}$,
  there exists $s_* \in [t - \max\{T, T_\varepsilon\}, t]$
  such that the two characteristics
  intersect, i.e.
  \begin{equation}\label{eq:117}
    X(s_*;t,x) = X(s_* ;t,x+\varepsilon ),
  \end{equation}
  which is absurd.

  The argument above, proved by contradiction, shows that $t \leq T$.
  Hence for all $t \in \reali_+$ and for all $x \in [0,1]$, the period
  of the orbit $\mathcal{O}\left(x, q(t,x)\right)$ must be greater
  than $t$.

  Since the period of the orbit $\mathcal{O}\left(x, q(t,x)\right)$
  diverges to $+\infty$ as $f(x, q(t,x)) \to \bar f$ for a suitable
  $\bar f \in \mathcal{F}$, letting $t \to + \infty$, we
  have~\eqref{eq:87}, proving \underline{Step~6}.

  \paragraph{\underline{Step~7}: $q$ asymptotically-in-time
    almost uniformly converges to $q_\infty$.}
  Recall assumption~\eqref{eq:138} and observe that the cases
  \begin{displaymath}
    r_o \equiv \underline{q}_R
    \quad \text{or} \quad
    r_o \equiv \overline{q}_R
  \end{displaymath}
  are trivial: since $\underline{q}_R$ (resp. $\overline{q}_R$)
  is a stationary solution
  to the equation in~\eqref{eq:75},
  we have that $q(t,x) = \underline{q}_R(x)$
  (resp. $q(t,x) = \overline{q}_R(x)$)
  for all $t \in \reali_+$ and $x \in [0,1]$.
  Using notation~\eqref{eq:33},
  this proves that $q$ (trivially)
  converges to the traveling wave
  $r_\infty$, with $A=\Su$
  (resp. $A=\emptyset$).

  Hence, we suppose
  \begin{equation}\label{eq:89}
    r_o \not\equiv \underline{q}_R
    \quad \text{and} \quad
    r_o \not\equiv \overline{q}_R\,.
  \end{equation}
  We consider two cases, depending on the number of points of global
  maximum of $\Sigma {\rho}_o$.

  \paragraph{\underline{Case~7.A}: $\Sigma{\rho}_o$ admits a unique
    point of global maximum in $\Su$.}
  Using notation~\eqref{eq:36}, we have
  $n=1$, $a_0=b_0=0$
  and $a_1=b_1=1$.

  By the characterization~\eqref{eq:106} and~\ref{eq:101}, this implies
  \begin{equation}
    \label{eq:152}
    \underline{q}_R(0)=\underline{q}_R(1) = \overline{q}_R(0)= \overline{q}_R(1)
    \quad\text{and}\quad
    \underline{q}_R(x) \neq \overline{q}_R(x)
    \ \ \mbox{for all}\ \ x \in \mathopen]0, 1 \mathclose[\,.
  \end{equation}
  Moreover, since~\eqref{eq:89} and~\eqref{eq:152} hold, there exists
  unique $\bar \xi \in \mathopen]0,1[$ such that
  \begin{equation}\label{eq:92bis}
    \int_{\Su} \left[
      r_o(x) - \underline{q}_R(x)
      \right] \d{x}
    =
    \int_{\bar \xi}^1 \left[
      \overline{q}_R(x) - \underline{q}_R(x)
      \right] \d{x}\,.
  \end{equation}

  \medskip
  \noindent
  Set $I \coloneq \{1,\ldots,2 \, \sharp(\mathcal F)\}$. Since
  $\sharp(\mathcal F) < +\infty$, the analysis of the level sets of
  $f$ developed in \textbf{\underline{Step~5}} yields that, for every
  $x \in [0,1]$, there exist at most $2\, \sharp(\mathcal F)$ points
  $q_i^x \in \mathcal D$, $i \in I$, such that
  $f(x, q_i^x) \in \mathcal F$; set
  $\mathcal Q^x := \{q_i^x : i \in I\}$.  For $q \in \reali$ we write
  $\operatorname{dist}(q,\mathcal Q^x) := \min_{i \in I} \modulo{q -
    q_i^x}$.

  By \textbf{\underline{Step~6}}, we have
  \begin{equation}\label{eq:88}
    \lim_{t \to + \infty}
    \sup_{x \in [0,1]}
    \operatorname{dist}\big(q(t,x), \mathcal Q^x\big) = 0\,.
  \end{equation}

  \begin{figure}[ht]
    \centering
    \begin{minipage}{0.5\textwidth}
      \centering
      \includegraphics[width=\linewidth]{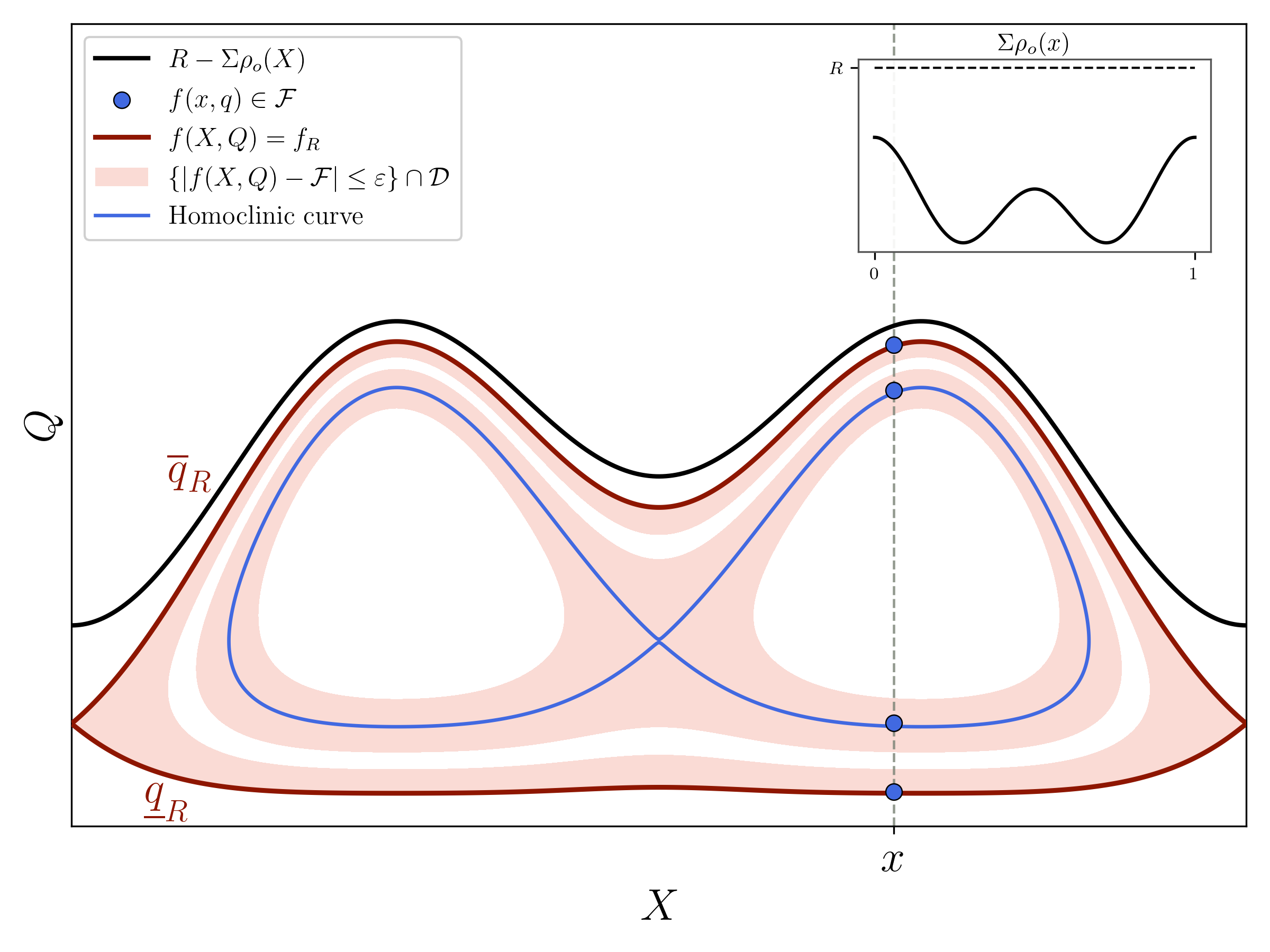}
    \end{minipage}%
    \hfill
    \begin{minipage}{0.45\textwidth}
      \caption{In the case $\Sigma\widehat{\rho}_o$ as in upper image,
        $\sharp({\mathcal F}) = 2$. For $x \in [0,1]$ fixed, we plot
        $(x, q_i^x)$ for $i \in \{1,4\}$ such that
        $f(x, q_i^x) \in \mathcal{F}$ (blue dots).  We highlight the
        region
        $\{(X,Q) \colon \min_{\bar f \in \mathcal{F}} \modulo{f(X,Q)
          {- }\bar f} < \varepsilon \} \cap \mathcal{D}$ for a
        positive $\epsilon$.  }
      \label{fig:5}
    \end{minipage}
  \end{figure}
  Fix $\varepsilon > 0$.  Using~\eqref{eq:87} and~\eqref{eq:88}, let
  $\bar T = \bar T(\varepsilon) >0$ be such that
  \begin{equation}\label{eq:80}
    \min_{\bar f \in \mathcal{F}}
    \modulo{f(x, q(t,x)) - \bar f} \leq \varepsilon
    \quad \text{and} \quad
    \operatorname{dist}\big(q(t,x), \mathcal Q^x\big) \leq \varepsilon
    \quad
    \text{for all } x \in [0, 1],\ t \geq \bar T
  \end{equation}
  and introduce
  \begin{displaymath}
    \mathcal{D}_\varepsilon \coloneqq
    \left\{ (X, Q) \in \mathcal{D}:
    \min_{\bar f \in \mathcal{F}}
    \modulo{f(X, Q) - \bar f} \leq \varepsilon
    \quad
    \text{and} \quad
    \operatorname{dist}\big(Q, \mathcal Q^X\big) \leq \varepsilon
    \right\}\,,
  \end{displaymath}
  so that~\eqref{eq:80} implies that
  $q(t,x) \in \mathcal{D}_\varepsilon$ for all $x \in [0,1]$ and
  $t \geq \bar T$.  Assume that $\varepsilon$ is so small that the
  sets
  $f^{-1}\big(\mathopen]\bar f - \varepsilon, \bar f +
  \varepsilon\mathclose[\big)$, for $\bar f \in \mathcal F$, are
  pairwise disjoint subsets of $\mathcal D$ and that
  \begin{eqnarray}
    \label{eq:108}
    \epsilon
    & <
    & \overline{q}_R(\bar \xi) - \underline{q}_R(\bar \xi)
      \quad
      \mbox{where } \bar \xi \mbox{
      satisfies~\eqref{eq:92bis} ;
      }
    \\
    \label{eq:83}
    \varepsilon
    & <
    & \frac{1}{4}
      \left(\min_{\bar f \in \mathcal{F} \setminus\{f_R\}} \bar f - f_R\right)
      \quad
      \mbox{if }
      \sharp(\mathcal{F}) \geq 2;
    \\
    \label{eq:90}
    \mathfrak{L}(\mathcal{D}_\varepsilon)
    & <
    & \min \left(
      \int_{\Su} \left(
      r_o(x) - \underline{q}_R(x)
      \right) \d{x},
      \int_{\Su} \left(
      \overline{q}_R(x) - r_o(x)
      \right) \d{x}
      \right),
  \end{eqnarray}
  $\mathfrak{L}$ being the Lebesgue measure. Recalling that
  $f_R = \min_{\bar f \in \mathcal{F}} \bar f$,
  by~\eqref{eq:158}--\eqref{eq:107}, and that~\eqref{eq:138}
  and~\eqref{eq:89} hold, the right hand sides of the previous
  inequalities are strictly positive.

  Fix $t \geq \bar T$. We record two preliminary observations:
  \begin{enumerate}[label=\textup{\textbf{($\mathcal{O}$\arabic*)}},
      ref=\textup{\textbf{($\mathcal{O}$\arabic*)}}]
    \item\label{item:7} There exists maximal $\delta >0$ such that for
          all
          $x \in [0, \delta \mathclose[ \cup \mathopen] 1 - \delta, 1 ]$,
          \begin{equation}
            \label{eq:81}
            \modulo{f\left(x, q(t,x)\right) - f_R}\leq \varepsilon
            \,,\quad
            \modulo{q(t, x) - \underline{q}_R(x)}\leq \varepsilon
            \mbox{ and }
            \modulo{q(t, x) - \overline{q}_R(x)}\leq \varepsilon\,.
          \end{equation}
          In fact, let
          $f_R^+ \coloneqq \min_{\bar f \in \mathcal{F} \setminus\{f_R\}}
            \bar f$, so that the curve $f^{-1}(f_R^+)$ is the maximal
          homoclinic/heteroclinic curve contained in the region
          $f^{-1}(\mathopen]f_R, +\infty \mathclose[)$.  Since it holds
          that
          $\{0,1\} = \mbox{argmax}_{[0,1]} \Sigma \widehat{\rho}_o$, we have
          that the set $f^{-1}(\mathopen]f_R, f_R^+[)$ is exactly the region
          confined by the curves $f^{-1}(f_R^+)$ and $f^{-1}(f_R)$ and it is
          homeomorphic to an annulus.  This proves~\eqref{eq:81}, first row.
          The continuity of $\underline{q}_R$ and $\overline{q}_R$
          and~\eqref{eq:152} ensure the second row.
    \item\label{item:8} Let $\bar x \in [\delta,1-\delta]$.  Then, the
          following holds:
            \begin{eqnarray}
              \label{eq:84}
              \modulo{q(t, \bar x ) - \underline{q}_R(\bar x)} \leq
              \varepsilon
              & \implies
              & \modulo{q(t, y) - \underline{q}_R(y)} \leq \varepsilon
              \quad \mbox{ for all }
              y \in [0, \bar x]\, ;
              \\
              \label{eq:85}
              \modulo{q(t, \bar x ) - \overline{q}_R(\bar x)} \leq
              \varepsilon
              & \implies
              & \modulo{q(t, y ) - \overline{q}_R(y)} \leq \varepsilon
              \quad \mbox{ for all }
              y \in [\bar x, 1].
            \end{eqnarray}
          This immediately follows from~\eqref{eq:80} and the fact that
          $x \mapsto q(t,x)$ can only feature positive jumps.
  \end{enumerate}
  \noindent Consider the set $\mathcal{J}_t$ of the points of jump
  discontinuity $x$ of $q(t)$ for which the traces $q(t, x-)$ and
  $q(t, x+)$ are close to different values of $\mathcal{Q}^x$, namely
  \begin{displaymath}
    \mathcal{J}_t
    \coloneqq
    \left\{ x \in [\delta,1-\delta] \colon
    q(t, x - ) < q(t, x +)
    \mbox{ and }
    \begin{array}{@{}r@{}r}
      \modulo{q(t, x-) - q_l^x} & < \varepsilon
      \\
      \modulo{q(t, x+) - q_r^x} & < \varepsilon
    \end{array}
    \;
    q_l^x, q_r^x \in \mathcal{Q}^x, \, q_l^x \neq q_r^x
    \right\}.
  \end{displaymath}
  Since $q(t) \in \BV([0, 1]; \reali)$,~\eqref{eq:80} holds
  and $|\mathcal{Q}^x| < + \infty$ for all $x \in [0,1]$,
  then $|\mathcal{J}_t| < + \infty$.

  Moreover, $\mathcal{J}_t \neq \emptyset$: otherwise,
  combining~\eqref{eq:80} and~\eqref{eq:81} we have
  \begin{displaymath}
    \begin{array}{lcl}
      \mbox{either }
      \sup_{x \in [0,1]} \modulo{q(t,x) - \underline{q}_R(x)} \leq \varepsilon
      & \mbox{ or }
      & \sup_{x \in [0,1]}
        \modulo{q(t,x) - \overline{q}_R(x)} \leq \varepsilon
      \\
      \mbox{either }
      \int_{\Su} \left[
      q(t,x) - \underline{q}_R(x)
      \right] \d{x} \leq \mathfrak{L}(\mathcal{D}_\varepsilon)
      & \mbox{ or }
      & \int_{\Su} \left[
        \overline{q}_R(x) - q(t,x)
        \right] \d{x} \leq \mathfrak{L}(\mathcal{D}_\varepsilon)
      \\
      \mbox{either }
      \int_{\Su} \left[
      r_o(x) - \underline{q}_R(x)
      \right] \d{x} \leq \mathfrak{L}(\mathcal{D}_\varepsilon)
      & \mbox{ or }
      & \int_{\Su} \left[
        \overline{q}_R(x) - r_o(x)
        \right] \d{x} \leq \mathfrak{L}(\mathcal{D}_\varepsilon)\,.
    \end{array}
  \end{displaymath}
    where we used the conservation of the mass on $\Su$. The latter
    statement contradicts~\eqref{eq:90}.

    \medskip We claim that there exists $\bar T_* \geq \bar T$ such
    that for all $t \geq \bar T_*$ there exists a unique
    $\bar \xi_t \in [\delta,1-\delta]$ belonging to $\mathcal{J}_t$.
    Moreover,
  \begin{equation}\label{eq:82}
    \modulo{q(t, x) - \underline{q}_R(x)} {\leq} \varepsilon
    \; \mbox{for all} \; x \in [0, \bar \xi_t \mathclose[
    \quad \mbox{and} \quad
    \modulo{q(t, x) - \overline{q}_R(x)} {\leq} \varepsilon
    \; \mbox{for all} \; x \in [\bar \xi_t, 1]\,.
  \end{equation}
  Using~\ref{item:7} and~\ref{item:8}, consider the unique non trivial
  case: there exists $x \in \mathcal{J}_t$ such that
  $q^x_l \neq \underline{q}_R(x)$ or $q^x_r \neq \overline{q}_R(x)$.

  Since $q(t)$ can allow only positive jumps, there exists
  $0 \leq \zeta_o \leq x \leq \theta_o \leq 1$ such that
  \begin{gather*}
    \modulo{q(t,x) - \underline{q}_R(x)} \leq \varepsilon \; \mbox{for all } x \in [0, \zeta_o \mathclose[\, ,
    \quad
    \modulo{q(t,x) - \overline{q}_R(x)} \leq \varepsilon \; \mbox{for all } x \in [\theta_o, 1]
    \\
    q(t,x) \in f^{-1} \left(
      \bigcup_{\bar f \in \mathcal{F} \setminus\{f_R\}}
      \mathopen] \bar f - \varepsilon, \bar f + \varepsilon \mathclose[
    \right)
    \mbox{ for all } x \in [\zeta_o, \theta_o\mathclose[.
  \end{gather*}
  Observe that $\zeta_o, \theta_o \in \mathcal{J}_t$.

  By~\cite[Theorem~11.9.4]{dafermos}, there exist unique forward
  generalized characteristics
  $[s \mapsto \zeta(s), s  \geq t]$,
  respectively $[s \mapsto \theta(s), s  \geq t]$,
  emanating from the point $(t, \zeta_o)$, respectively
  $(t, \theta_o)$.
  Moreover, for all $s \geq t$,
  \begin{equation}\label{eq:161}
    q(s, \zeta(s)-) < q(s, \zeta(s)+)
    \quad
    \mbox{and}
    \quad
    q(s, \theta(s)-) < q(s, \theta(s)+)
  \end{equation}
  and the speed of the characteristics is given by the Rankine-Hugoniot
  condition
  \begin{eqnarray*}
    \dot \zeta(s)
    & =
    & \frac{f\left(\zeta(s), q(s,\zeta(s)-)\right)
      -
      f\left( \zeta(s), q(s,\zeta(s)+)\right)
      }{
      q(s, \zeta(s)- )- q(s, \zeta(s)+)
      }
    \\
    \dot \theta(s)
    & =
    & \frac{f\left(\theta(s), q(s,\theta(s)-)\right)
      -
      f\left( \theta(s), q(s,\theta(s)+)\right)
      }{
      q(s, \theta(s)- )- q(s, \theta(s)+)
      }\,.
  \end{eqnarray*}

  \medskip
  \noindent
  We claim that as long as $\zeta$ and $\theta$ do not
  collide with other shocks, there exists
  $\alpha >0$ such that
  $\dot \zeta(s) > \alpha$ and $\dot \theta(s) < -\alpha$
  for all $s \geq t$.

  In fact,
  recall that
  \begin{equation}\label{eq:86}
    \modulo{f(\zeta_o,q(t,\zeta_o -)) - f_R} \leq \varepsilon,
    \quad \mbox{and} \quad
    \modulo{f(\zeta_o,q(t,\zeta_o +)) - f_\zeta} \leq \varepsilon,
  \end{equation}
  for suitable $f_\zeta \in \mathcal{F} \setminus \{f_R\}$.
  Hence, we have
  \begin{align*}
    \dot \zeta(t)
     &
    =
    \frac{
      f\left( \zeta_o, q(t,\zeta_o+)\right)
      - f\left(\zeta_o, q(t,\zeta_o-)\right)
    }{
      q(t, \zeta_o+) - q(t, \zeta_o - )
    }
    \\
     &
    \geq
    \frac{
      f_{\zeta}  - \varepsilon - f_R - \varepsilon
    }
    {q(t, \zeta_o+) - q(t, \zeta_o - )
    }
    \tag*{[By~\eqref{eq:86}]}
    \\
     &
    >
    \underbrace{\frac{
        \min_{\bar f \in \mathcal{F} \setminus\{f_R\}} \bar f
        - f_R
      }
      {2\,R}\,.
    }_{\alpha}
    \tag*{[By~\eqref{eq:83} and $ q(t, \zeta_o+) - q(t, \zeta_o - ) \leq R$ ]}
  \end{align*}
  As long as $\zeta$ does not collide with another shock,
  by~\cite[Theorem~11.9.6]{dafermos}, for all $s \geq t$ each point
  $(s, \zeta(s))$ is a point of continuity of the map
  $[(\tau, x) \mapsto q(\tau, x-)]$ on the region
  $\{ (\tau, x) : 0 < \tau < + \infty, x \leq \zeta(\tau)\}$ and also
  a point of continuity of the map $[(\tau, x) \mapsto q(\tau, x+)]$
  on the region
  $\{ (\tau, x) : 0 < \tau < + \infty, x \geq \zeta(\tau)\}$.

  Due to the continuity of $f$,
  the maps
  $[s \mapsto f(\zeta(s), q(s,\zeta(s)\pm ))]$ are
  also
  continuous.

  Hence, the maps $[s \mapsto f(\zeta(s), q(s,\zeta(s)\pm ))]$
  cannot pass from being
  $\varepsilon$-close to one value of $\mathcal F$ to being $\varepsilon$-close
  to a different one without violating~\eqref{eq:80} at some intermediate
  instant.

  Altogether, we have
  \begin{displaymath}
    \modulo{ f\left(\zeta(s), q(s,\zeta(s)-)\right) - f_R} \leq \varepsilon
    \quad \mbox{and} \quad
    \modulo{ f\left(\zeta(s), q(s,\zeta(s)+)\right) - f_\zeta} \leq \varepsilon
  \end{displaymath}
  and, consequently, $\dot \zeta(s) > \alpha > 0$
  for all $s \geq t$ as long as $\zeta$ does not collide with another shock.

  Analogously, we have that $\dot \theta(s) < -\alpha < 0$,
  for all $s \geq t$ as long as $\theta$ does not collide with another shock.

  \medskip
  Since $\dot \zeta (s) > \alpha$ and $\dot \theta (s) < -\alpha$,
  either $ \zeta$
  and $\theta$ collide with one another, or one of them
  collides first with a shock
  emanated from a point in $\mathcal{J}_t$
  lying between $\zeta_o$ and $\theta_o$.
  In both cases there exists
  $t \leq t_1  < t + 1/(2\alpha) < +\infty$ such that
  one of the two cases occurs.

  After the collision, only one shock survives,
  $\sharp(\mathcal{J}_{t_1})$ has decreased by (at least) one and
  $q(t_1)$ is such that
  \begin{displaymath}
    \modulo{q(t_1,x) - \underline{q}_R(x)} \leq \varepsilon \quad \mbox{for all } x \in [0, \zeta(t_1) \mathclose[\, ,
    \quad
    \modulo{q(t_1,x) - \overline{q}_R(x)} \leq \varepsilon \quad \mbox{for all } x \in [\theta(t_1), 1]
  \end{displaymath}
  \begin{displaymath}
    q(t_1,x) \in f^{-1} \left(
      \bigcup_{\bar f \in \mathcal{F} \setminus\{f_R\}}
      \mathopen] \bar f - \varepsilon, \bar f + \varepsilon \mathclose[
    \right)
    \mbox{ for all } x \in [\zeta(t_1), \theta(t_1)\mathclose[.
  \end{displaymath}
  Repeting of the above argument a finite number of times proves the
  existence of $\bar T_* < +\infty$ and a unique shock curve
  $[ s \mapsto \xi(s)]$ for $s \geq \bar T_*$ such that
  \begin{displaymath}
    \modulo{q(s, x) - \underline{q}_R(x)} \leq \varepsilon \quad \mbox{for all } x \in [0, \xi(s) \mathclose[\, ,
    \quad
    \modulo{q(s, x) - \overline{q}_R(x)} \leq \varepsilon \quad \mbox{for all } x \in [\xi(s), 1]\,.
  \end{displaymath}

  \medskip
  \noindent
  The conservation of mass yields for all $s \geq \bar T_*$
  \begin{equation}\label{eq:94}
    \int_{\Su} \left[
      r_o(x) - \underline{q}_R(x)
      \right] \d{x}
    =
    \int_{\Su} \left[
      q(s,x) - \underline{q}_R(x)
      \right] \d{x} \,.
  \end{equation}
  Use~\eqref{eq:80} and~\eqref{eq:119} to get
  \begin{eqnarray*}
    \int_{\Su} \left[
    q(s,x) - \underline{q}_R(x)
    \right] \d{x}
    & =
    & \int_{0}^{\xi(s)} \left[
      q(s,x) - \underline{q}_R(x)
      \right] \d{x}
      +
      \int_{\xi(s)}^1 \left[
      q(s,x) - \underline{q}_R(x)
      \right] \d{x}
    \\
    & \leq
    & \varepsilon \, \xi(s) + \int_{\xi(s)}^1 \left[
      \overline{q}_R(x) - \underline{q}_R(x)
      \right] \d{x}
    \\
    & \leq
    & \varepsilon  + \int_{\xi(s)}^1 \left[
      \overline{q}_R(x) - \underline{q}_R(x)
      \right] \d{x}\,, \qquad \mbox{ and}
    \\
    \int_{\Su} \left[
    q(s,x) - \underline{q}_R(x)
    \right] \d{x}
    & \geq
    & \int_{\xi(s)}^1 \left[
      \overline{q}_R(x) - \underline{q}_R(x)
      \right] \d{x} - \varepsilon (1 - \xi(s))
    \\
    & \geq
    & \int_{\xi(s)}^1 \left[
      \overline{q}_R(x) - \underline{q}_R(x)
      \right] \d{x}
      - \varepsilon\, ,
  \end{eqnarray*}
  so that
  \begin{equation}\label{eq:91}
    \modulo{
      \int_{\Su} \left[
        q(s,x) - \underline{q}_R(x)
        \right] \d{x}
      -
      \int_{\xi(s)}^1 \left[
        \overline{q}_R(x) - \underline{q}_R(x)
        \right] \d{x}
    } \leq \varepsilon \,.
  \end{equation}
  Recall $\bar \xi$ from~\eqref{eq:92bis}.
  Combine~\eqref{eq:94}$,$~\eqref{eq:91} and~\eqref{eq:92bis} to obtain
  \begin{displaymath}
    \modulo{
      \int_{\bar \xi}^{\xi(s)} \left[
        \overline{q}_R(x) - \underline{q}_R(x)
        \right] \d{x}
    } \leq \varepsilon \,
  \end{displaymath}
  Using~\eqref{eq:108} and the fact that the map
  $\left[y \mapsto \int_{\bar \xi}^y \left( \overline{q}_R(x) -
      \underline{q}_R(x) \right) \d{x}\right]$ for $y \in [0,1]$ is
  strictly increasing, the arbitrariness of $\varepsilon$ implies that
  $\xi(s) \underset{s \to + \infty}{\longrightarrow} \bar \xi$.
  \medskip Altogether, we have that $q(t)$ asymptotically-in-time
  almost uniformly converges to $q_\infty$ where
  \begin{displaymath}
    q_\infty(x)
    =
    \begin{cases}
      \underline{q}_R(x)
      & x \in [0, \bar \xi[
      \\
      \overline{q}_R(x)
      & x \in [\bar \xi, 1]
    \end{cases}
  \end{displaymath}
  and $\bar \xi$ is the unique point in $\mathopen]0,1[$
  satisfying~\eqref{eq:92bis}\,.

  \paragraph{\underline{Case~7.B}: The set of points of
    global maxima of $\Sigma \rho_o$ is a finite union of closed intervals.}
  The characterization~\eqref{eq:106}, together with~\eqref{eq:119},
  implies that
  \begin{equation}\label{eq:23}
    \underline{q}_R(x) = q(t, x) = \overline{q}_R(x)
    \quad
    \mbox{for all
    $x \in \mathop{\rm argmax}_{[0,1]} \Sigma \widehat{\rho}_o$,
    $t \in \reali_+$}.
  \end{equation}
  Hence, for all $i \in \{0, \ldots, n\}$, the curves
  $[s\in\reali_+ \to a_i]$ and $[s\in\reali_+ \to b_i]$ are maximal
  divides, in the sense of~\cite[Definition~10.3.3]{dafermos}. They
  divide the domain $\reali_+\times [0,1]$ in $2n+1$ disjoint regions
  $\mathcal{D}_i$ and $\mathcal{C}_i$
  \begin{displaymath}
    \mathcal{C}_i = \mathopen]0,+\infty\mathclose[ \times \mathopen]a_i, b_i\mathclose[\,,
      \quad
      \mathcal{D}_i = \mathopen]0,+\infty\mathclose[ \times \mathopen]b_i, a_{i+1}\mathclose[\,.
  \end{displaymath}

  For all $i \in \{0, \ldots, n-1\}$ and
  $(t,x) \in \mathcal{D}_i$,
  the associated maximal backward generalized
  characteristic $[s \mapsto X(s, x, q(t,x)), s \in [0,t]]$
  cannot intersect the maximal divides $[s\in\reali_+ \to a_i]$ and
  $[s\in\reali_+ \to b_i]$
  and lies in $\mathcal{D}_i$ during
  its whole evolution.
  So,
  the argument developed in \textbf{\underline{Case~7.A}}
  holds in each region
  $\mathcal{D}_i$:
  there exists unique
  $\bar \xi_i \in [b_i, a_{i+1}]$
  and
  \begin{equation}
    \label{eq:156bis}
    q_\infty(x) = \begin{cases}
      \underline{q}_R(x)
       &
      x \in [b_i, \bar \xi_i \mathclose[ \\
      \overline{q}_R(x)
       & x \in [\bar \xi_i, a_{i+1}]\,
    \end{cases}
    \quad \mbox{for all $x \in [b_i, a_{i+1}]$}
  \end{equation}
  such that $q(t)$ asymptotically-in-time almost uniformly
  converges to $q_\infty$ on $[b_i, a_{i+1}]$.
  Moreover, since $\underline{q}_R(x) < \overline{q}_R(x)$
  for all $x \in ]b_i, a_{i+1}[$,
  $\bar{\xi}_i$ is uniquely determined by
  \begin{displaymath}
    \int_{b_i}^{\bar \xi_i} \underline{q}_R(x)\d{x}
    + \int_{\bar \xi_i}^{a_{i+1}} \overline{q}_R(x)\d{x}
    =
    \int_{b_i}^{a_{i+1}} r_o(x) \d{x}\,.
  \end{displaymath}
  On the other hand, $q(t)$ is stationary on $[a_i, b_i]$
  by~\eqref{eq:23}, so it trivially
  asymptotically-in-time converges to $\underline{q}_R = \overline{q}_R$
  in that interval.

  Altogether, $q$ asymptotically-in-time almost
  uniformly converges to $q_\infty$ as in~\eqref{eq:33bis}.

  \paragraph{\underline{Step~8}: $r$ asymptotically-in-time
    almost uniformly converges to $r_\infty$.}
  Fix any sequence $t_\nu \in \reali_+ $ such that
  $\lim_{\nu \to +\infty} t_\nu = + \infty$.
  By \textbf{\underline{Case~7.B}}, $q$,
  the solution to~\eqref{eq:75}, asymptotically-in-time almost uniformly
  converges to $q_\infty$ as in~\eqref{eq:33bis},
  i.e.
  for all $\varepsilon >0$
  \begin{displaymath}
    \lim_{\nu \to +\infty}
    \sup_{
      \{
      x \in [0,1] \colon
      \underset{i \in \{0, \ldots, n-1\}}{\min}
      \modulo{x-\bar \xi_i} \geq \varepsilon
      \}
    }
    \modulo{q(t_\nu, x) - q_\infty(x)} =0\,.
  \end{displaymath}
  Since $q$ is periodic in the $x$-variable, we have
  \begin{align*}
    \!\!\!\!\!
    \lim_{\nu \to +\infty}
    \sup_{
      \{x \in \reali
      \colon
      \underset{k \in \interi,  i \in \{0, \ldots, n-1\}}{\min}
      \modulo{x-\bar \xi_i-k} \geq \varepsilon\}
    }
    \modulo{q(t_\nu, x) - \widehat q_\infty(x)}
     &
    =0
    \\
    \!\!\!\!\!
    \lim_{\nu \to +\infty}
    \sup_{
      \{x \in \reali
      \colon
      \underset{k \in \interi,  i \in \{0, \ldots, n-1\}}{\min}
      \modulo{x-\bar \xi_i-k} \geq \varepsilon\}
    }
    \modulo{\widehat{r}(t_\nu, x{-}\lambda t_\nu) {- }\widehat q_\infty(x)}
     &
    =0
    \tag*{[\mbox{$q(t_\nu, x) = \widehat{r}(t_\nu, x+\lambda t_\nu)$}]}
    \\
    \!\!\!\!\!
    \lim_{\nu \to +\infty}
    \sup_{
      \{\xi \in \reali
      \colon
      \underset{k \in \interi,  i \in \{0, \ldots, n-1\}}{\min}
      \modulo{\xi {- }\lambda t_\nu {-}\bar \xi_i{-}k} \geq \varepsilon\}
    }
    \modulo{\widehat{r}(t_\nu, \xi) {- }\widehat{q}_\infty(\xi {-} \lambda t_\nu)}
     &
    =0
    \tag*{[\mbox{Set $\xi = x + \lambda t_\nu$}]},
  \end{align*}
  concluding the proof.

  \paragraph{\underline{Step~9}: $(\rho, r_\infty)$
    solves~\eqref{eq:3} with initial datum $(\rho_o, q_\infty)$ in the
    sense of \Cref{def:sol}.}
  Observe that
  $(\rho, r_\infty) \in \C0\left(\reali_+; \W11(\Su; \reali^k) \times
    \L1(\Su; \reali)\right)$ due to the continuity of translations in
  $\W11(\Su; \reali^k) \times \L1(\Su; \reali)$.

  With reference to the notation in Definition~\ref{def:sol} define
  \begin{displaymath}
    w (t,x,u)
    \coloneqq
     v_L\left(\Sigma\widehat{\rho (t)} \, (x) + u \right)\,,
      \quad \mbox{ so that } \quad
    w (t,x,u) =
    v_L\left(\Sigma\widehat{\rho_o} (x-\lambda t )+ u \right)\,.
  \end{displaymath}
  Let $\kappa \in \reali$ ,
  $\psi \in \Cc1(\mathopen]0, +\infty \mathclose[ \times \reali;
  \reali_+)$ and suppose that
  $\spt(\psi) \subset [t_1, t_2] \times [x_1, x_2]$ for some
  $0 < t_1 < t_2 < +\infty$, $-\infty < x_1 < x_2 < +\infty$.
  \begin{align}
    \label{eq:31}
     &
    \int_{t_1}^{t_2} \int_{x_1}^{x_2} \modulo{\widehat{r}_\infty(t,x) - \kappa} \, \partial_t \psi(t,x)
    \d{x}\d{t}
    \\
    \nonumber
     &
    +
    \int_{t_1}^{t_2}\int_{x_1}^{x_2} \sgn(\widehat{r}_\infty(t,x) - \kappa) \,
    \left[
      \widehat{r}_\infty(t,x) \, w\left(t,x,\widehat{r}_\infty (t,x)\right)
      - \kappa \, w(t,x,\kappa)
      \right] \partial_x \psi(t,x) \d{x}\d{t}
    \\
    \nonumber
     &
    -
    \int_{t_1}^{t_2} \int_{x_1}^{x_2} \sgn(\widehat{r}_\infty(t,x) - \kappa) \,
    \kappa \, \partial_x w(t,x,\kappa) \, \psi(t,x) \d{x}\d{t}
    \\
    \nonumber
    =
     &
    \int_{t_1}^{t_2} \int_{x_1-\lambda t }^{x_2 - \lambda t }
    \modulo{\widehat{q}_\infty(\xi) - \kappa}
    \left(
    \partial_t \Psi(t, \xi) - \lambda \partial_\xi \Psi(t, \xi)
    \right)
    \d{\xi}\d{t}
    \\
    \nonumber
     &
    +
    \int_{t_1}^{t_2}\!\!\!
    \int_{x_1-\lambda t }^{x_2 - \lambda t }\!
    \sgn(\widehat{q}_\infty(\xi) {-} \kappa) \,
    \left[
      \widehat{q}_\infty(\xi) \, v_L\left(\Sigma\widehat{\rho_o} (\xi ){+ }\widehat{q}_\infty(\xi) \right)
      {- }\kappa \, v_L\left(\Sigma\widehat{\rho_o} (\xi){+} \kappa \right)
      \right]
    \partial_\xi \Psi(t, \xi) \d{\xi}\d{t}
    \\
    \nonumber
     &
    -
    \int_{t_1}^{t_2} \int_{x_1-\lambda t }^{x_2 - \lambda t }
    \sgn(\widehat{q}_\infty(\xi) - \kappa) \,
    \kappa \, v_L'\left(\Sigma\widehat{\rho_o} (\xi )+ \kappa \right) \Sigma \widehat{\rho}_o'(\xi)
    \Psi(t,\xi) \d{\xi}\d{t}
  \end{align}
  where we have used the change of variable $\xi = x-\lambda t$
  and set $\Psi(t, \xi) \coloneqq \psi(t, \xi + \lambda t)$.
  Recall the definition of $f$~\eqref{eq:75} in \textbf{\underline{Step~1}},
  so that
  \begin{eqnarray*}
    \widehat{q}_\infty(\xi) \, v_L\left(\Sigma\widehat{\rho_o} (\xi ){+ }\widehat{q}_\infty(\xi) \right)
    & =
    & f(\xi, \widehat{q}_\infty(\xi)) + \lambda\, \widehat{q}_\infty(\xi)
    \\
    \kappa \, v_L\left(\Sigma\widehat{\rho_o} (\xi){+} \kappa \right)
    & =
    & f(\xi, \kappa) + \lambda \kappa
    \\
    \kappa \, v_L'\left(\Sigma\widehat{\rho_o} (\xi )+ \kappa \right) \Sigma \widehat{\rho}_o'(\xi)
    & =
    & \partial_X f(\xi, \kappa)\,.
  \end{eqnarray*}
  Altogether, we have
  \begin{align}
    \nonumber
    \mbox{\eqref{eq:31}}
    =
     &
    \int_{t_1}^{t_2} \int_{x_1-\lambda t }^{x_2 - \lambda t }
    \underbrace{
      \modulo{\widehat{q}_\infty(\xi) - \kappa}
    }_{A_\infty(\xi)}
    \partial_t \Psi(t, \xi)
    \d{\xi}\d{t}
    \\
    \nonumber
     &
    +
    \int_{t_1}^{t_2}
    \int_{x_1-\lambda t }^{x_2 - \lambda t }
    \underbrace{
      \sgn(\widehat{q}_\infty(\xi) {-} \kappa) \,
      \left[
        f(\xi, \widehat{q}_\infty(\xi))
        -
        f(\xi, \kappa)
        \right]
    }_{B_\infty(\xi)}
    \partial_\xi \Psi(t, \xi) \d{\xi}\d{t}
    \\
    \nonumber
     &
    -
    \int_{t_1}^{t_2} \int_{x_1-\lambda t }^{x_2 - \lambda t }
    \underbrace{
      \sgn(\widehat{q}_\infty(\xi) - \kappa) \,
      \partial_X f(\xi, \kappa)
    }_{C_\infty(\xi)}
    \Psi(t,\xi) \d{\xi}\d{t}.
  \end{align}
  Since $(\rho,r)$ solves~\eqref{eq:3} in the sense of~\Cref{def:sol},
  for all
  $\psi_*\in \Cc1(\mathopen]0, +\infty \mathclose[ \times \reali;
  \reali_+)$,
  \begin{align}
    \nonumber
     &
    \int_0^{+\infty} \int_{\reali} \modulo{\widehat{r}(t,x) - \kappa} \, \partial_t \psi_*(t,x)
    \d{x}\d{t}
    \\
    \nonumber
     &
    +
    \int_0^{+\infty}\int_{\reali} \sgn(\widehat{r}(t,x) - \kappa) \,
    \left[
      \widehat{r}(t,x) \, w\left(t,x,\widehat{r} (t,x)\right)
      - \kappa \, w(t,x,\kappa)
      \right] \partial_x \psi_*(t,x) \d{x}\d{t}
    \\
    \nonumber
     &
    -
    \int_0^{+\infty} \int_{\reali} \sgn(\widehat{r}(t,x) - \kappa) \,
    \kappa \, \partial_x w(t,x,\kappa) \,\psi_*(t,x) \d{x}\d{t}
    \\
    \nonumber
    =
     &
    \int_0^{+\infty} \int_{\reali}
    \underbrace{
      \modulo{\widehat{r}(t,\xi+\lambda t ) - \kappa}
    }_{A(t,\xi)}
    \partial_t \Psi_*(t, \xi)
    \d{\xi}\d{t}
    \\
    \nonumber
     &
    +
    \int_0^{+\infty}
    \int_{\reali}
    \underbrace{
      \sgn(\widehat{r}(t,\xi+\lambda t ) {-} \kappa) \,
      \left[
        f(\xi, \widehat{r}(t,\xi+\lambda t ))
        -
        f(\xi, \kappa)
        \right]
    }_{B(t,\xi)}
    \partial_\xi \Psi_*(t, \xi) \d{\xi}\d{t}
    \\
    \nonumber
     &
    -
    \int_0^{+\infty} \int_{\reali}
    \underbrace{
      \sgn(\widehat{r}(t,\xi+\lambda t ) - \kappa) \,
      \partial_X f(\xi, \kappa)
    }_{C(t,\xi)}
    \Psi_*(t,\xi) \d{\xi}\d{t}
    \\
    \label{eq:34}
    \geq
     & 0
  \end{align}
  where again we have used $\xi \coloneqq x - \lambda t$ and
  $\Psi_*(t,\xi) \coloneqq \psi_*(t,\xi+\lambda t) $.

  \textbf{\underline{Step~7}} ensures that
  for a.e. $\xi \in \reali$,
  $\widehat{r}(t,\xi+\lambda t ) = q(t,\xi )
    \underset{t \to + \infty}{\longrightarrow} \widehat{q}_\infty(\xi)$
  almost uniformly.
  Since the map $[u \mapsto \modulo{u-\kappa}]$ and
  $[u \mapsto \sgn(u-\kappa)\left(
        f(\xi, u) - f(\xi, \kappa)
        \right)]$ are continuous, we have that
  for a.e. $\xi \in \reali$,
  \begin{displaymath}
    A(t)\underset{t \to + \infty}{\longrightarrow} A_\infty
    \quad
    B(t)\underset{t \to + \infty}{\longrightarrow} B_\infty
    \quad \mbox{a.e. and in } \Lloc1(\reali; \reali).
  \end{displaymath}
  We claim that also
  \begin{displaymath}
    C(t)\underset{t \to + \infty}{\longrightarrow} C_\infty
    \quad \mbox{a.e. and in } \Lloc1(\reali; \reali).
  \end{displaymath}
  We tackle the unique non trivial case
  $E_\kappa \coloneqq \{ \xi \in \reali \colon \widehat{q}_\infty(\xi) = \kappa\} \neq \emptyset$.
  For all $\xi \in E_\kappa$ we obtain by~\ref{item:43}
  in \Cref{thm:Claudia}
  \begin{align*}
    f(\xi, \kappa)
    =
    \widehat{q}_\infty(\xi)
    \left(
    v_L\left(\Sigma\widehat{\rho_o} (\xi )+ \widehat{q}_\infty(\xi) \right)
    - \lambda
    \right)
    =
    f_R\,.
  \end{align*}
  which implies that
  $\partial_X f(\xi, \kappa) =0$ for a.e. $\xi \in \reali$
  with $\widehat{q}_\infty(\xi) = \kappa$
  by~\cite[Theorem~4 in Section~4.2]{MR4907811}.
  Hence, for a.e. $\xi \in \reali$ with $\widehat{q}_\infty(\xi) = \kappa$
  we have that
  \begin{displaymath}
    C(t, \xi) =  \sgn(\widehat{r}(t,\xi+\lambda t ) - \kappa) \,
    \partial_X f(\xi, \kappa)
    =
    0
    =
    C_\infty(\xi),
  \end{displaymath}
  yielding the claim.

  Fix $\varepsilon >0$. Hence, there exists $T \in \reali_+$ such that
  for all $t \in \reali$ with $t \geq T$,
  \begin{equation}\label{eq:29bis}
    \norma{A(t) - A_\infty}_{\L1(K; \reali)} +
    \norma{B(t) - B_\infty}_{\L1(K; \reali)}
    +\norma{C(t) - C_\infty}_{\L1(K; \reali)}
    \leq \varepsilon
  \end{equation}
  where
  $K \coloneqq [x_1 - \modulo{\lambda } t_2, x_2 + \modulo{\lambda}
  t_2]$.  Therefore, we have
  \begin{eqnarray*}
    &
    & \left|
      \int_{t_1}^{t_2} \int_{x_1-\lambda t }^{x_2 - \lambda t }
      \left[
      {A_\infty(\xi)}
      \partial_t \Psi(t,\xi)
      +
      {B_\infty(\xi)}
      \partial_\xi \Psi(t,\xi )
      -
      {C_\infty(\xi)}
      \Psi(t,\xi )
      \right] \d{\xi}\d{t}
      \right.
    \\
    &
    & \quad {- } \left.
      \int_{t_1}^{t_2}\!\!\! \int_{x_1-\lambda t }^{x_2 - \lambda t }\!
      \left[
      {A(t{+}T,\xi)}
      \partial_t \Psi(t,\xi)
      {+}
      {B(t{+}T,\xi)}
      \partial_\xi \Psi(t,\xi )
      {-}
      {C(t{+}T,\xi)}
      \Psi(t,\xi )
      \right] \d{\xi}\d{t}
      \right|
    \\
    & \leq
    & \modulo{
      \int_{t_1+T}^{t_2+T} \int_K
      \left(
      {A_\infty(\xi)} - A(t, \xi)
      \right)
      \partial_t \Psi(t-T,\xi)
      \d{\xi}\d{t}
      }
    \\
    &
    & \quad +
      \modulo{
      \int_{t_1+T}^{t_2+T} \int_K
      \left(
      {B_\infty(\xi)} - B(t, \xi)
      \right)
      \partial_\xi \Psi(t-T, \xi)
      \d{\xi}\d{t}
      }
    \\
    &
    & \quad
      +
      \modulo{
      \int_{t_1+T}^{t_2+T} \int_K
      \left(
      {C_\infty(\xi)} - C(t, \xi)
      \right)
      \Psi(t-T, \xi)
      \d{\xi}\d{t}
      }
    \\
    & \leq
    &
      \left(
      \norma{\partial_t \Psi}_{\L\infty(\reali_+ \times \reali; \reali)}
      +
      \norma{\partial_\xi \Psi}_{\L\infty(\reali_+ \times \reali; \reali)}
      +
      \norma{\Psi}_{\L\infty(\reali_+ \times \reali; \reali)}
      \right)(t_2-t_1)\, \varepsilon
  \end{eqnarray*}

  The arbitrariness of $\varepsilon$ and~\eqref{eq:34} (with
  $\Psi_*(t, \xi) = \Psi(t-T, \xi)$)
  imply that~\eqref{eq:31} is nonnegative, showing the entropy
  inequality in \Cref{def:sol}.

  The verification that $\rho$ satisfies the equations in \Cref{def:sol}
  with $r_\infty$ in place of $r$ is a straightforward adaptation of the
  above argument, using that each $\rho^i$ is a $\lambda$-traveling wave
  and that $v_{NL}^i \in \C0$.

  This completes the proof of \textbf{\underline{Step~9}} and of
  Theorem~\ref{thm:Claudia}.
\end{proofof}

\appendix

\section{Functional Analysis Results}
\label{sec:technique}

\begin{lemma}
  \label{lem:CBV}
  Let $G \in \BVloc (\reali; \reali)$. Then, there exists a function
  $H \in \C0 (\reali, \reali)$ such that
  \begin{displaymath}
    \mbox{for a.e. } x \in \reali
    \qquad
    \sgn H (x) = \sgn G (x) \,.
  \end{displaymath}
\end{lemma}

\begin{proofof}{Lemma~\ref{lem:CBV}}
  $G$ is $\BVloc$ and right continuous, hence it is continuous on
  $\reali$ except at countably many points. Define $\mathcal{G}^+$ as
  the interior of $\left\{ x \in \reali \colon G (x) >0 \right\}$ and
  $\mathcal{G}^-$ as the interior of
  $\left\{ x \in \reali \colon G (x) <0 \right\}$. For every open
  subset $A$ of $\reali$, introduce the continuous map
  $\Upsilon_A (x) \coloneq \min\left\{ 1, d (x, \reali\setminus A)
  \right\}$. Then, one can choose
  $H (x) \coloneq \Upsilon_{\mathcal{G}^+} (x) -
  \Upsilon_{\mathcal{G}^-} (x)$. Indeed, if $x$ is a point of
  continuity of $G$ and $G (x) >0$, then $x \in \mathcal{G}^+$ and
  $H (x) >0$. A similar observation applies if $x$ is a point of
  continuity of $G$ and $G (x) <0$. If $x$ is a point of continuity of
  $G$ and $G (x) =0$, then
  $x \in \reali \setminus (\mathcal{G}^+ \cup \mathcal{G}^-)$, so that
  $H (x) = 0$.
\end{proofof}

\begin{lemma}
  \label{lem:sign}
  Let $F \in \W11(\reali; \reali)$ and
  $G \in \BVloc(\reali; \reali)$ such that for all $x \in \reali$
  \begin{equation}\label{eq:131}
    G(x) = 0
    \quad
    \implies
    \quad
    F(x) = 0\,.
  \end{equation}
  Then, the map
  $\left[ x \mapsto \sgn\left(G(x)\right)\, F(x) \right]$ is in
  $\W11 (\reali; \reali)$ and
  \begin{equation}
    \label{eq:129}
    \dfrac{\d~}{\d{x}} \left( \sgn\left(G(x)\right)\, F(x) \right)
    =
    \sgn\left(G(x)\right) \, \dfrac{\d~}{\d{x}} F(x)\,.
  \end{equation}
\end{lemma}

\begin{proofof}{Lemma~\ref{lem:sign}}
  We have to prove that for all $\phi \in \Cc1(\reali; \reali)$,
  \begin{align*}
    \int_\reali \sgn\left(G(x)\right) \, F(x) \, \phi'(x)
    \d{x}
    =
    - \int_\reali \sgn\left(G(x)\right) \, F'(x) \, \phi(x)
    \d{x} \,.
  \end{align*}
  By Lemma~\ref{lem:CBV}, we can replace $G$ above with a
  $H \in \C0 (\reali; \reali)$ such that $\sgn H (x) = \sgn G (x)$ for
  a.e.~$x \in \reali$.  Denote $\mathring{\spt(\phi)}$ the
  interior of ${\spt(\phi)}$ and compute:
  \begin{eqnarray}
    \label{eq:136}
    &
    & \int_\reali \sgn\left(H(x)\right)
      \left(F(x) \, \phi'(x) + F'(x) \, \phi(x)\right)
      \d{x}
    \\
    \nonumber
    & =
    & \int_\reali \sgn\left(H(x)\right)
      \left(
      F(x) \,
      \phi(x)
      \right)'
      \d{x}
    \\
    \nonumber
    & =
    & \int_{\left\{x \in \reali \colon H(x)>0\right\} \, \cap \, \mathring{\spt(\phi)}}
      \left(
      F(x) \,
      \phi(x)
      \right)'
      \d{x}
      -
      \int_{\left\{x \in \reali \colon H(x)<0 \right\}\, \cap \, \mathring{\spt(\phi)}}
      \left(
      F(x) \,
      \phi(x)
      \right)'
      \d{x}\,.
  \end{eqnarray}
  Since $H \in \C0(\reali; \reali)$, the set
  $\left\{ x \in \reali \colon H(x) > 0 \right\}$ ia a countable union
  of disjoint open intervals and $H$ vanishes at each of the extreme
  points of these intervals. Similarly, $\pint{\spt\phi}$ is the
  countable union of disjoint open intervals and $\phi$ vanishes at
  each of the extreme points of these intervals. Hence, by~\eqref{eq:131}
  \begin{displaymath}
    \left\{ x \in \reali \colon H(x) >0 \right\}
    \, \cap \, \mathring{\spt(\phi)}
    = \bigcup_{i \in \naturali}
    \,
    \mathopen]a_i, b_i\mathclose[
    \quad \mbox{ with } \quad
    F (a_i) \, \phi (a_i) = 0
    \mbox{ and }
    F (b_i) \, \phi (b_i) = 0 \,.
  \end{displaymath}
  Clearly, see for instance~\cite[Corollary~8.10]{Brezis},
  $F \, \phi \in \W11(\reali; \reali)$ so that
  \begin{eqnarray*}
    \int_{\left\{ x \in \reali \colon H(x)>0 \right\}
    \, \cap \, \mathring{\spt(\phi)}
    }
    \left(
    F(x) \,
    \phi(x)
    \right)'
    \d{x}
    & =
    & \sum_{i \in \naturali}
      \int_{a_i}^{b_i}
      \left(
      F(x) \,
      \phi(x)
      \right)'
      \d{x}
    \\
    & =
    & \sum_{i \in \naturali}
      \left[F(b_i)\phi(b_i) - F(a_i)\phi(a_i) \right] = 0
  \end{eqnarray*}
  The same argument applies to
  $\left\{ x \in \reali \colon H(x) <0 \right\}\, \cap \,
  \mathring{\spt(\phi)}$.  Hence,~\eqref{eq:136} vanishes. It is
  then immediate to verify that
  $\left[ x \mapsto \sgn \left( H(x)\right) F'(x) \right]$ is in
  $\L1\left(\reali; \reali\right)$, concluding the proof.
\end{proofof}

\begin{lemma}
  \label{lem:variational}
  Fix a non trivial open interval $I$ and let
  $f \in \Lloc1(I; \reali)$ such that
  \begin{equation}
    \label{eq:73}
    \int_I f(t) \, \chi(t) \d{t} = 0
    \quad
    \mbox{for all}
    \quad
    \chi \in \Cc1(I; \reali)
    \mbox{ with }
    \int_I \chi(t) \d{t} = 1 \, .
  \end{equation}
  Then $f=0$ a.e.~in $I$.
\end{lemma}

\begin{proofof}{Lemma~\ref{lem:variational}}
  At first, we claim that~\eqref{eq:73} holds also for all
  $\chi \in \Cc0(I; \reali)$ such that $\int_I \chi(t) \d{t} = 1$.  In
  fact, introduce the sequence of standard mollifiers
  \begin{displaymath}
    \rho_n \in \Cc\infty(\reali; \reali)
    \quad \mbox{ such that }
    \spt(\rho_n) \subset \left[-\frac{1}{n}, \frac{1}{n}\right],
    \quad
    \int_\reali \rho_n(t) \d{t} = 1 ,
    \quad
    \rho_n \geq 0.
  \end{displaymath}
  For $n \in \naturali$ sufficiently large,
  $\chi_n \coloneqq \chi* \rho_n \in \Cc1(I; \reali)$,
  $\int_I \chi_n(t) \d{t} = 1$ and
  $\chi_n \xrightarrow[n\to +\infty]{} \chi$ uniformly on compact sets
  (see~\cite[Proposition~4.21]{Brezis}). Hence, by~\eqref{eq:73} and
  the Dominated Convergence Theorem,
  $\int_I f(t) \bar{\chi}(t) \d{t} = \lim_{n\to+\infty} = \int_I f(t)
  \chi_n(t) \d{t} = 0$, proving the claim.

  By linearity, it is immediate to verify that~\eqref{eq:73} holds
  also for all $\chi \in \Cc0(I; \reali)$ such that
  $\int_I \chi(t) \d{t} \neq 0$. Fix now $\chi \in \Cc0(I; \reali)$
  such that $\int_I \chi(t) \d{t} = 0$. Introduce the sequence
  \begin{displaymath}
    \chi_n \in \Cc0(I; \reali)
    \quad \mbox{ by } \quad
    \chi_n (t) \coloneqq
    \begin{cases}
      \begin{array}{lcc}
        \bar{\chi}(t) \left(1 + \frac{1}{n}\right)
        & \mbox{if}
        & \bar{\chi}(t) >0
        \\
        \bar{\chi}(t)
        & \mbox{if}
        & \bar{\chi}(t) \leq 0 \,.
      \end{array}
    \end{cases}
  \end{displaymath}
  By construction $\chi_n \xrightarrow[n \to + \infty]{} \bar \chi$
  uniformly and $\int_I \chi_n(t) \d{t} \neq 0$ for all
  $n \in \naturali$. Again by the Dominated Convergence Theorem,
  $\int_I f(t) \, \bar{\chi}(t) \d{t} = \lim_{n\to+\infty} = \int_I
  f(t) \, \chi_n(t) \d{t} = 0$. We can now
  apply~\cite[Corollary~4.24]{Brezis} to conclude the proof.
\end{proofof}

\section*{Acknowledgments}
The authors were partly supported by the GNAMPA~2025 project
\emph{Modellli di Traffico, di Biologia e di Dinamica dei Gas Basati
  su Sistemi di Equazioni Iperboliche}. CN was supported by the PNRR
project \emph{Ricerca DM 118/2023}, CUP~F13C23000360007.

\section*{Conflict  of Interest}
The authors declare no conflicts of interest in this paper.

  {\small

    \bibliographystyle{abbrv}

    \bibliography{tw}

  }

\end{document}